\documentclass[12pt,letter]{scrartcl}

\usepackage{amsthm}
\usepackage{mathptmx}       % selects Times Roman as basic font
\usepackage{helvet}         % selects Helvetica as sans-serif font
\usepackage{courier}        % selects Courier as typewriter font
\usepackage{makeidx}         % allows index generation
\usepackage{multicol}        % used for the two-column index
\usepackage[bottom]{footmisc}% places footnotes at page bottom

\usepackage[pdftex]{graphicx}
\usepackage{float}
\usepackage{hyperref}
\usepackage{amssymb}
\usepackage{amsmath}

\newcommand{\be}{\begin{eqnarray}}
\newcommand{\ee}{\end{eqnarray}}
\newcommand{\bez}{\begin{eqnarray*}}
\newcommand{\eez}{\end{eqnarray*}}

\newcommand{\Hscalprod}[2]{{\langle #1,#2\rangle}}                         %scalar product
\newcommand{\Hset}[2]{\left\{#1\vphantom{#2}\right.\;\left|\;\vphantom{#1}#2\right\}}        %set

\newtheorem{theorem}{Theorem}[section]

\newtheorem{proposition}{Proposition}[section]

\theoremstyle{remark}
\newtheorem{remark}{Remark}[section]

\theoremstyle{definition}
\newtheorem{definition}{Definition}[section]
\newtheorem{example}{Example}[section]
\newtheorem{problem}{Problem}[section]
\newtheorem*{acknowledgement}{Acknowledgement}

\makeindex             % used for the subject index
\begin{document}

\date{October 2011}

\title{Permutahedra and Associahedra}
\subtitle{Generalized associahedra from the geometry of finite reflection groups}
\author{Christophe Hohlweg\thanks{D\'epartement de math\'ematiques, 
Universit\'e du Qu\'ebec \'a Montr\'eal, 
Case postale 8888, succ. Centre-ville
Montr\'eal (Qu\'ebec), H3C 3P8
Canada, 
 \url{hohlweg.christophe@uqam.ca}}}

\maketitle

\abstract{This is a chapter in an upcoming Tamari Festscrift. Permutahedra are a class of convex polytopes arising naturally from the study of finite reflection groups, while generalized associahedra are a class of polytopes indexed by finite reflection groups. We present the intimate links those two classes of polytopes share.
%\\
%Keywords: Coxeter groups; permutahedron ; Tamari lattice; reflection groups ; Cambrian fans; Cambrian lattices; %Generalized associahedron; cluster. \\
%MSC 20F55 ; 06B99 ; 52B11 ; 05E99
}

\section{Introduction}
\label{H_sec:Introduction}

The purpose of this survey is to explain the realization of generalized associahedra and Cambrian lattices (which are generalizations of the Tamari lattice) from the geometrical point of view of finite reflection groups. 

The story of the associahedron starts in 1963 when J.~Stasheff~\cite{H_St63}, while studying homotopy theory of loop spaces, constructed a cell complex whose vertices correspond to the possible compositions of $n$ binary operations. This cell complex turns out to be the boundary complex of a convex polytope: {\em the associahedron}. Long forgotten was that D.~Tamari considered in his 1951 thesis a lattice which graph is the graph of the associahedron, and for which he had a realization in dimension $3$  (see in this regard J.-L.~Loday's text in this volume~\cite[\S8]{H_Lo11}).

One of the easiest ways to realize an associahedron is to cut out a standard simplex. It turns out that by cutting out this associahedron, we obtain the {\em classical permutahedron}. The permutahedron is a polytope that arises from the symmetric group seen as a finite reflection group. This construction owed to S.~Shnider and S.~Sternberg~\cite{H_ShSt94} (see also~\cite{H_St95}), and later completed by J.-L.~Loday~\cite{H_Lo04}, builds a bridge between the classical permutahedron and the associahedron. It carries many combinatorial and geometrical properties: for instance, this transformation maps the weak order on the symmetric group to the Tamari lattice.

A similar object, {\em the cyclohedron}, was later discovered by R.~Bott and C.~Taubes in 1994~\cite{H_BoTa94} in connection with knot theory. Realizations of the cyclohedron were given by M.~Markl~\cite{H_Ma99}, R.~Simion~\cite{H_Si03} and V.~Reiner~\cite{H_Rei02}, but none of these realizations exhibit a link with the symmetric group or other finite reflection groups.  In 2003, S.~Fomin and A.~Zelevinsky~\cite{H_FoZe03} discovered, in their study of finite type cluster algebras, a family of polytopes, indexed by finite reflection groups, that contains the associahedron and the cyclohedron.  These {\em generalized associahedra} were first realized by F.~Chapoton, S.~Fomin and A.~Zelevinsky~\cite{H_ChFoZe03} but still, these realizations were not obtained from a permutahedron of the corresponding finite reflection group.

In 2007, C.~Hohlweg and C.~Lange~\cite{H_HoLa07}, and subsequently with H.~Thomas~\cite{H_HoLaTh11}, described many realizations of {\em generalized associahedra}, that are obtained by `removing some facets' from the permutahedron of the corresponding finite reflection group (see Figure~\ref{H_fig:move} below): we start from a finite reflection group $W$, a special ordering of the simple reflections and a permutahedron for $W$ and end with a realization of a generalized associahedra. This way of realizing generalized associahedra has many benefits: it maps the weak order on $W$ to Cambrian lattices, the vertices of these realizations are labeled by clusters, and their normal fans are Cambrian fans (from N.~Reading's Cambrian lattices and fans, see~\cite{H_Re11}).  The spinal cord of this construction is the Coxeter singletons that allow to pinpoint nicely the facets of the permutahedron that have to be removed.  

\begin{figure}
%\sidecaption
\begin{center}
\includegraphics[scale=0.4]{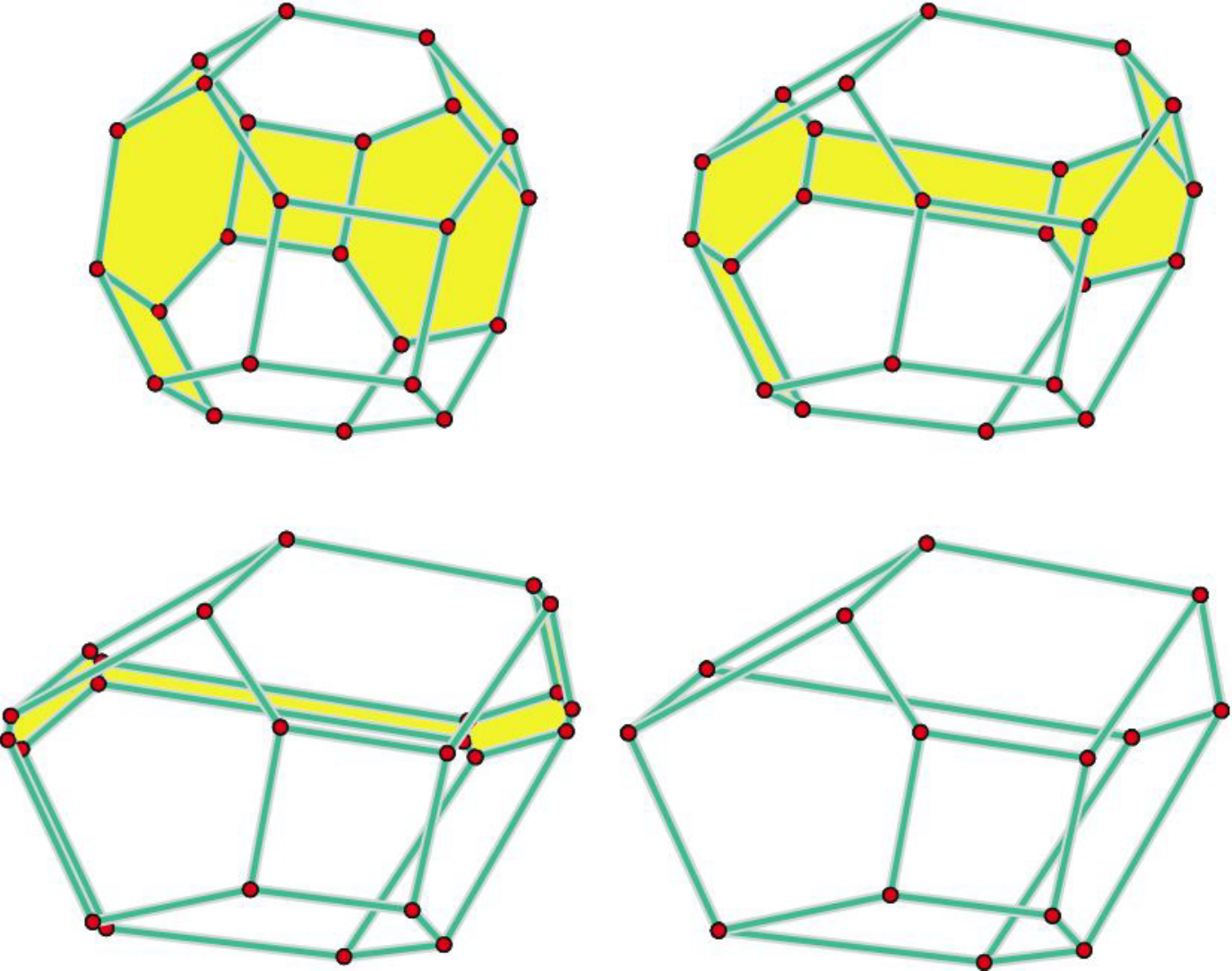}
\caption{The process of removing facets from a permutahedron in order to obtain a generalized associahedron, as shown in \cite{H_HoLa07}. The yellow facets represent the facets in the process of being removed.}
\label{H_fig:move}
\end{center}
\end{figure}

The first part of this survey is dedicated to permutahedra of finite reflection groups, and how they encode important data about the group. The second part is dedicated to realizations of generalized associahedra from a given permutahedron and how they result in a geometrical construction of Cambrian lattices and fans. Along the way, we give numerous examples and figures, and discuss open problems and further developments on the subject.

For more details on polytopes and fans, we refer to the book by G.~M.~Ziegler~\cite{H_Zi95}, from which we use the general notations. A nice presentation of finite reflection groups can be found in the  book of J.~Humphrey~\cite{H_Hu90}.

%%%%%%%%%%%%%%%%%%%%%%%%%
\section{Permutahedra and finite reflection groups}
%%%%%%%%%%%%%%%%%%%%%%%%%%

We consider a finite dimensional $\mathbb R$-Euclidean space~$(V,\Hscalprod{\cdot}{\cdot})$.

\subsection{Finite reflection groups} %%%%%%%%%
\index{reflection group (finite)}

 A {\em finite reflection group} is a finite subgroup of the orthogonal group $O(V)$ generated by reflections.  A {\em reflection}\footnote{We consider only orthogonal reflections in this text.} can be defined relative to the hyperplane it fixes pointwise, or by a normal vector to this hyperplane. Let us fix some notation: if $H$ is an hyperplane in $V$ and $\alpha\in V$ is a normal vector to $H$, the reflection $s_\alpha$ is the unique linear isometry which fixes $H=\alpha^\perp$ pointwise and maps $\alpha$ to $-\alpha$. A general formula for $s_\alpha(v)$, for a vector $v\in V$, follows:
\bez
  s_\alpha (v) = v- 2\frac{\Hscalprod{v}{\alpha}}{\Hscalprod{\alpha}{\alpha}}\alpha. 
\eez

\index{reflection group (finite)!dihedral group}

\begin{example} The basic example of finite reflection groups are dihedral groups (see Figure~\ref{H_fig:dihedral}). Take a regular $n$-gon $P$ in the affine plane $\mathbb R^2$ centered in $O=(0,0)$; the symmetry group of $P$ is the dihedral group $\mathcal D_n$, which is generated by the reflections it contains.  Each reflection of $\mathcal D_n$ is determined by an axis of symmetry of $P$; that is, a line passing through a vertex  $A$ and the point $O$, or the line passing through the middle of an edge and the point $O$. 
 \begin{figure}   %%%%%%%%%Figure
  \begin{center}
  \includegraphics[scale=0.6]{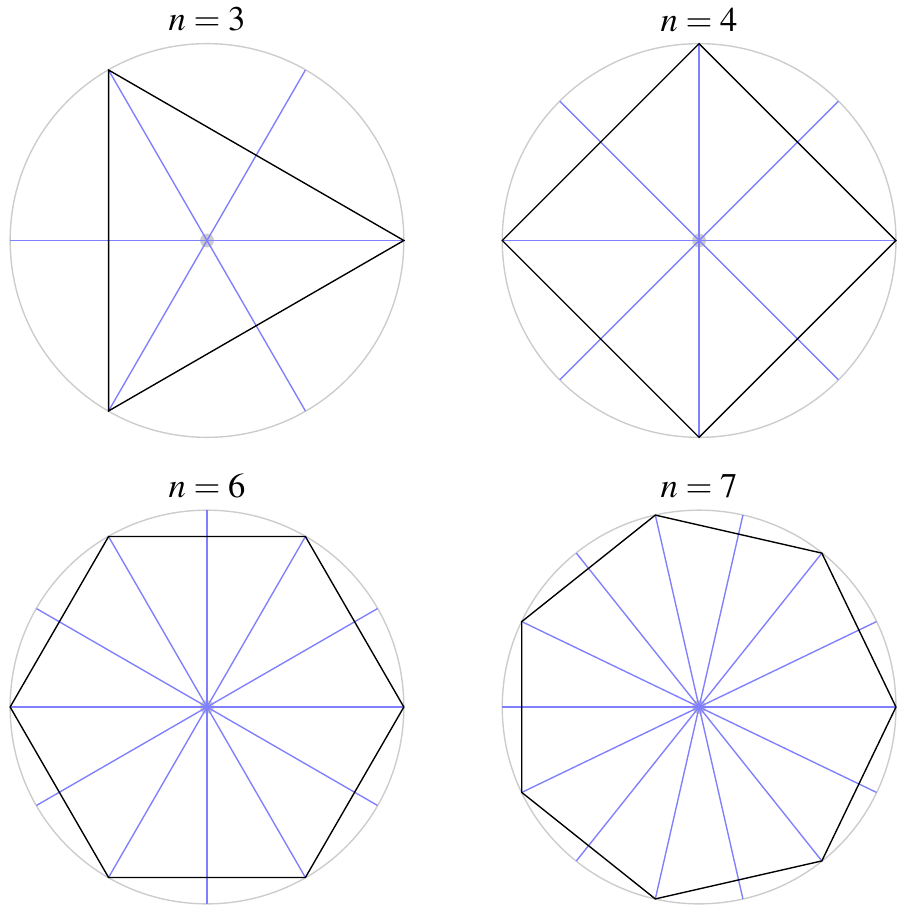}
    \caption{The triangle, the square, the hexagon and the heptagon with their axes of symmetry. Their isometry groups are reflection groups: the groups $\mathcal D_3$, $\mathcal D_4$, $\mathcal D_6$ and $\mathcal D_7$, respectively.}
  \label{H_fig:dihedral} 
  \end{center}
 \end{figure}   %%%%%%%%%%%%
\end{example}

\begin{example}
\label{H_ex:Sn}
\index{reflection group (finite)!symmetric group}
 Symmetric groups are certainly the most studied  finite reflection groups, since they enjoy a particular and very detailed combinatorial representation. They appear geometrically as isometry groups of standard regular simplexes. Let us  describe this representation. 
We consider $V=\mathbb R^n$, together with its canonical basis $\mathcal B=\{e_1,\dots,e_n \}$ on which the symmetric group $S_n$  acts linearly by permutation of the coordinates:
\bez
\sigma\big((x_1,\dots,x_n)\big) = (x_{\sigma(1)},\dots,x_{\sigma(n)}). 
\eez
The transposition $\tau_{ij}=(i\, j)$ ($1\leq i<j\leq n$) acts therefore as the orthogonal reflection that fixes pointwise the hyperplane $H_{i,j}=\{(x_1,\dots,x_n)\in\mathbb R^n\,|\, x_i =x_j\}$ and maps $\alpha_{i,j}:=e_{j}-e_i$ to $-\alpha_{i,j}=e_i-e_{j}$. Since $S_n$ is generated by transpositions, $S_n$ is a finite reflection group in $O(\mathbb R^n)$. It is well-known  that this representation of $S_n$, as a subgroup of $O(\mathbb R^n)$, is faithful, but not irreducible, since $S_n$ fixes the line spanned by the vector $(1,1,\dots,1)$.   
\end{example}

\begin{example}
\label{H_ex:Bn}
\index{reflection group (finite)!hyperoctahedral group}
Another instance of well-studied finite reflection groups are the {\em hyperoctahedral groups}. Hyperoctahedral groups also enjoy  a particular and very detailed combinatorial representation.  We consider $V=\mathbb R^n$, together with its canonical basis $\mathcal B=\{e_1,\dots,e_n \}$. Consider  the reflection $s_0:e_1\mapsto -e_1$ and recall the notations of Example~\ref{H_ex:Sn}. The group $W'_n$ generated by the reflection $s_0$ together with the reflections $\tau_{i,j}$ is a finite reflection group called {\em a hyperoctahedral group}, which contains $S_n$ as a reflection subgroup. For more on the special combinatorics of hyperoctahedral groups, see the book by A.~Bj\"orner and F.~Brenti  \cite{H_BjBr05}.
\end{example}

\index{reflection group (finite)!symmetry group of the dodecahedron}
\begin{example} A last and interesting example is the isometry group of the dodecahedron. The {\em regular dodecahedron}  is a $3$-dimensional convex polytope composed of $12$ regular pentagonal faces.  Its isometry group is denoted\footnote{The notation $W(H_3)$ refers to the classification of finite reflection groups, see Remark~\ref{H_re:classi}.} by $W(H_3)$. This group is a reflection group spanned by the reflections associated to the planes of symmetry of the dodecahedron.
\end{example}

\subsection{Permutahedra as $\mathcal V$-polytopes} %%%%%%%%
\label{H_ss:perm}

A natural way to study a given reflection group $W$ in $O(V)$ is to consider the $W$-orbit of a point $\pmb a\in V$. Here $V$ is endowed with its natural structure of affine Euclidean space. To ensure we get all the information we want, we need to choose this point generically: $\pmb a\in V$ is {\em generic} if it is not fixed by any reflections in $W$, or equivalently, if $\pmb a$ is not a point in a hyperplane corresponding to a reflection of $W$. Such a choice of a point $\pmb a$ is always possible since $W$ is finite. 

The idea behind {\em permutahedra} is to study $W$ with tools from polytope theory. Instead of considering only the $W$-orbit of $\pmb a$, we consider the polytope obtained as the convex hull of this orbit. 

\index{permutahedron}
\begin{definition} Let $\pmb a$ be a generic point of $V$.  The {\em Permutahedron} $\textnormal{Perm}^{\pmb a} (W)$ is the $\mathcal V$-polytope\footnote{The polytope is given as the convex hull of a set of points.} obtained as the convex hull of the $W$-orbit of $\pmb a$:
\bez
\textnormal{Perm}^{\pmb a}(W)= \textnormal{conv}\,\Hset{w(\pmb a)}{w\in W}.
\eez
This class of polytopes is called {\em Coxeter permutahedra} and is sometimes referred to by $W$-permutahedra in the literature. As we will see, the combinatorics of $\textnormal{Perm}^{\pmb a}(W)$ does not depend of the choice of $\pmb a$, as long as this point is generic.
\end{definition}

Since $W$ is finite, the $W$-orbit of $\pmb a$ is finite and $\textnormal{Perm}^{\pmb a}(W)$ is a polytope with a rich structure of faces. Before studying  $\textnormal{Perm}^{\pmb a}(W)$ in more detail, let us cover some examples. 

 \index{permutahedron! of dimension $2$}
\begin{example} 
\label{H_ex:PermD4}
Permutahedra of dimension $2$ arise from dihedral groups\footnote{The dimension of a polytope is the dimension of the affine space it spans.}: a permutahedron for $\mathcal D_n$ is a $2n$-gon. 

For instance, take $W=\mathcal D_4$, the symmetry group of a square (see Figure~\ref{H_fig:permdihedral}). The axes of symmetries, which correspond to the reflections in $W$, are the diagonals of the square and the lines passing through the middle of two opposite edges.  We pick a point $\pmb a$ that is not located in these lines. The red points represent the $W$-orbit of  $\pmb a$, and $\textnormal{Perm}^{\pmb a}(W)$ is the pale red convex octagon. Observe that the number of vertices of $\textnormal{Perm}^{\pmb a}(W)$ is equal to $8=|W|$, the cardinality of $W$. 
 \begin{figure}  %%%%%%%%%Figure
 \begin{center}
   \includegraphics{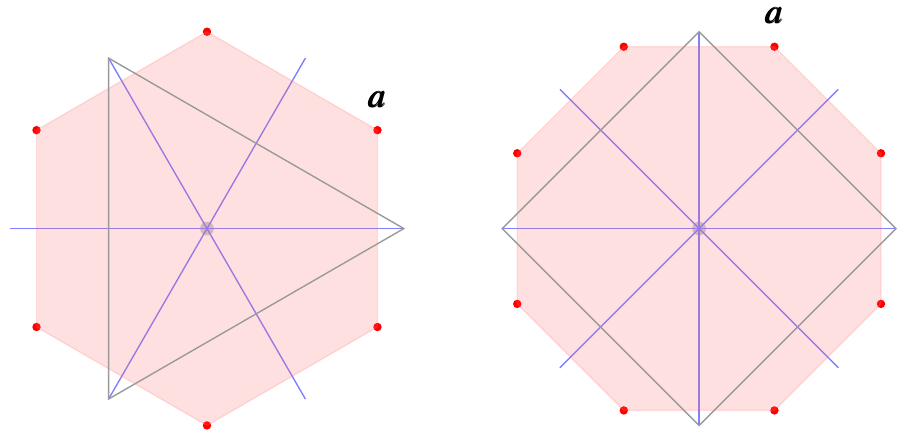}
   \caption{The polygon on the left is a permutahedron  of the dihedral group $\mathcal D_3$ and the polygon on the right is a permutahedron  of the dihedral group $\mathcal D_4$}
 \label{H_fig:permdihedral} 
 \end{center}
 \end{figure}   %%%%%%%%%
 \end{example}

 \begin{example}
 \label{H_ex:SnPerm}
 \index{permutahedron!classical permutahedron}
 {\bf The classical permutahedron}\\

The well studied {\em classical permutahedron} is provided in the framework of symmetric groups (see for instance \cite{H_Zi95,H_HoLa07} or \cite[\S2.2]{H_Lo04}).   The classical permutahedron\footnote{We use here the notation of G.~M.~Ziegler~\cite[Example 0.10]{H_Zi95}.} $\Pi_{n-1}$ is defined as the convex hull of all the permutations of the vector
$\pmb a=(1,2,\dots, n)\in\mathbb R^n$:
\bez
\Pi_{n-1}= \textnormal{conv}\,\{(\sigma(1),\dots,\sigma(n))\in\mathbb R^n\,|\,\sigma \in S_n\}.
\eez
It is a $n-1$-dimensional simple\footnote{A polytope is $d$-dimensional simple if any vertex is contained in precisely $d$ facets. This property is very strong: the face lattice of a simple polytope is completely determined by looking at its vertices and edges (see for instance \cite[\S3.4]{H_Zi95}).} convex polytope which lives in the affine hyperplane:
\bez
V_{\pmb a} = \Hset{(x_1,\dots,x_n)\in\mathbb R^n}{\sum_{i=1}^n x_i=\frac{n(n+1)}{2}}.
\eez
The vertices are naturally labeled by permutations of $S_n$ by denoting $M(\sigma)=(\sigma(1),\dots,\sigma(n))$.  The classical permutahedron $\Pi_2$ is shown in Figure~\ref{H_fig:a2_perm}. We observe that $\Pi_2$ is a permutahedron for the dihedral group $\mathcal D_3$, since $S_3$ acts as $\mathcal D_3$ on the affine hyperplane $V_{\pmb a}$.
 
 \begin{figure}
\begin{center}
\includegraphics[scale=0.5]{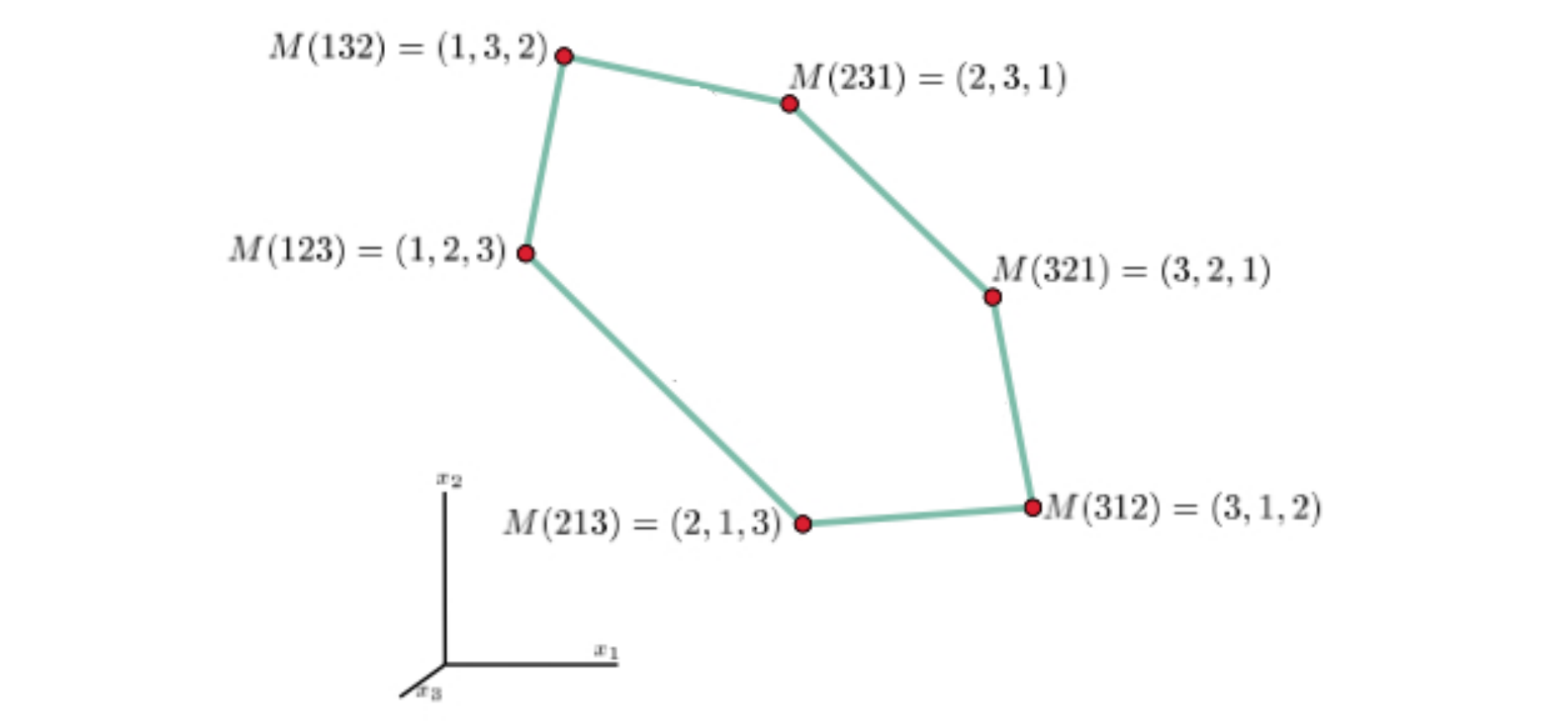}
\caption{The classical permutahedron $\Pi_2$ as shown in \cite{H_HoLa07}. This $2$-dimensional permutahedron in $\mathbb R^3$ is shown in the affine hyperplane $V_{\pmb{a}}$ with $x_1+x_2+x_3=6$. }
\label{H_fig:a2_perm}
\end{center}
\end{figure}
\end{example}

 \index{permutahedron! of dimension $3$}
\begin{example} There are two kinds of permutahedra of dimension $3$. 
\label{H_ex:PermDim3}
The first kind is obtained in $\mathbb R^3$ from permutahedra of dimension $2$ by considering the isometry group of a regular polygonal prism: take a $n$-gon $P$ in $\mathbb R^2=\textnormal{span}\,\{e_1,e_2\}$; the convex hull of  $P$ and of $P+e_3$ is a {\em regular $n$-gonal prism}. Its isometry group $W$ is isomorphic to $\mathcal D_n \times \mathcal D_2$ and $\textnormal{Perm}^{\pmb a}(W)$ is a $2n$-gonal prism.

The second kind of permutahedra of dimension $3$ arises from the symmetric group $S_4$, the hyperoctahedral group $W'_3$ and the isometry group of the dodecahedron $W(H_3)$. Examples are shown in  Figure~\ref{H_fig:permDim3}.

\begin{figure}
\begin{center}
   \includegraphics[scale=0.7]{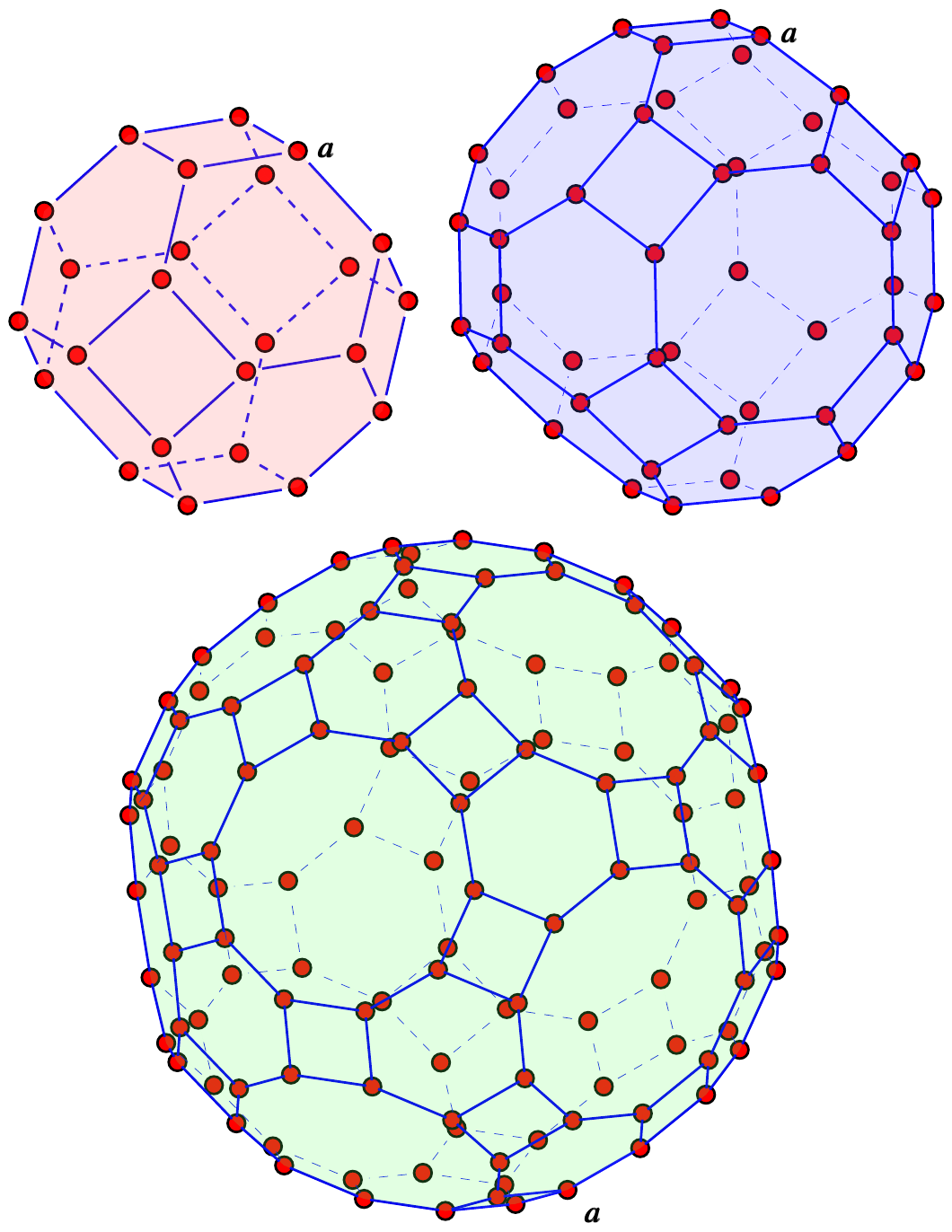}
 \caption{{\bf Permutahedra of dimension $\pmb 3$} From the top left to the bottom: a permutahedron for the symmetric group $S_4$, a permutahedron for the hyperoctahedral group $W'_3$ and a permutahedron for $W(H_3)$.}
  \label{H_fig:permDim3}
\end{center}
\end{figure}

\end{example}
 
  In both these examples, we can observe  that the number of vertices of the permutahedron is the cardinality of the finite reflection group $W$. This observation remains true in general.

\begin{theorem}
\label{H_thm:VertexPerm}
 The $W$-orbit of the generic point $\pmb a$ in $V$ is the vertex set of $\textnormal{Perm}^{\pmb a}(W)$.  In particular $|\textnormal{vert}\,(\textnormal{Perm}^{\pmb a}(W))|=|W|$.
\end{theorem}

 We observed in Figure~\ref{H_fig:permdihedral} that the reflection hyperplanes of a dihedral group $W=\mathcal D_m$, which are lines in this case,  cut the space into precisely $2m=|W|$ connected components which are polyhedral cones. In each of these cones lies a unique element of the $W$-orbit of $\pmb a$. Moreover, we observe that each edge of the permutahedron is directed by a normal vector to one of the hyperplanes of reflection.  A similar phenomenon can be observed for any finite reflection group. Looking at normal vectors of the reflecting hyperplanes of $W$ is the key to a thorough study of permutahedra, which will imply the statement in Theorem~\ref{H_thm:VertexPerm}.

\subsection{Root systems and permutahedra as $\mathcal H$-polytopes}  %%%%%%

We just described permutahedra as a $\mathcal V$-polytope, that is, as the convex hull of a set of points. Another way to describe a polytope is as an $\mathcal H$-polytope, that is, as the intersection of half-spaces. In order to do so, we introduce another very important tool related to finite Coxeter groups: {\em root systems}.

\subsubsection{Root systems}

We consider again a finite reflection group $W$ acting on $V$. As explained before, a reflection $s$ is uniquely determined by a given hyperplane $H$, or by a given normal vector $\alpha$ of $H$ and we write $s=s_\alpha$. This second point of view brings us to consider finite subsets $\Phi$ of $V$  such that 
\bez
  W=\langle s_\alpha\,|\, \alpha \in \Phi \rangle .
\eez 

\begin{definition}
\label{H_def:root}
\index{roots!root system}
 A {\em root system} for $W$ is a finite and nonempty subset $\Phi$ of nonzero vectors of $V$ that enjoys the following property:
\begin{enumerate}
\item $W=\langle s_\alpha\,|\,\alpha\in \Phi\rangle$;
\item for any $\alpha\in \Phi$, the line $\mathbb R\alpha$ intersects $\Phi$ in only the two vectors $-\alpha$ and $\alpha$;
\item $\Phi$ is stable under the action of $W$. 
\end{enumerate}
\index{roots}
The elements of $\Phi$ are called {\em roots}.
\end{definition}

\begin{remark}
\begin{enumerate}
\item Roots are normal vectors for reflection hyperplanes associated to $W$.
\item Root systems  exist for all finite reflection groups. Indeed, the set
\bez
\Phi:=\Hset{\pm \alpha}{s_\alpha\in W\textrm{ and } \Hscalprod{\alpha}{\alpha}=1}
\eez
obviously verifies the two first properties of Definition~\ref{H_def:root}. Moreover, it is not difficult to check that if $\alpha,\beta$ are two nonzero vectors of $V$ then 
$s_\alpha s_\beta s_\alpha= s_{s_\alpha(\beta)}.$
Therefore, $\Phi$ is stable under the action of $W$.  However it is important to note that roots do not have to have the same length in general: see \cite{H_Hu90} for more details. 
\end{enumerate}  
\end{remark}

\begin{example}
\label{H_ex:rootD4} 
Consider $\mathbb R^2$ with its canonical basis $e_1,e_2$. For the dihedral group $\mathcal D_3$, there are three axes of symmetries which are directed by vectors $e_1$, $e_1+\sqrt 3 e_2$ and $-e_1+\sqrt 3e_2$. Let us set $\alpha_1=\sqrt 3e_1+e_2$ and $\alpha_2=-\sqrt 3e_1+e_2$,  then a root system  for $\mathcal D_3$ is
\bez
\Phi=\{\pm \alpha_1,\pm \alpha_2,\pm(\alpha_1+\alpha_2)\}.
\eez
Note that the roots all have the same length (see Figure~\ref{H_fig:rootdihedral}).

Consider now the dihedral group $\mathcal D_4$. Since the axes of symmetries are pairwise orthogonal, they are all directed by normal vectors corresponding to the reflections of $\mathcal D_4$. Therefore, a root system  for $\mathcal D_4$ is
\bez
\Phi=\{\pm e_1,\pm e_2,\pm(e_1+e_2), \pm(e_2-e_1)\}.
\eez
Note that the roots do not all have the same length (illustrated by two colors in Figure~\ref{H_fig:rootdihedral}).

 \begin{figure}  %%%%%%%%%Figure
 \begin{center}
   \includegraphics{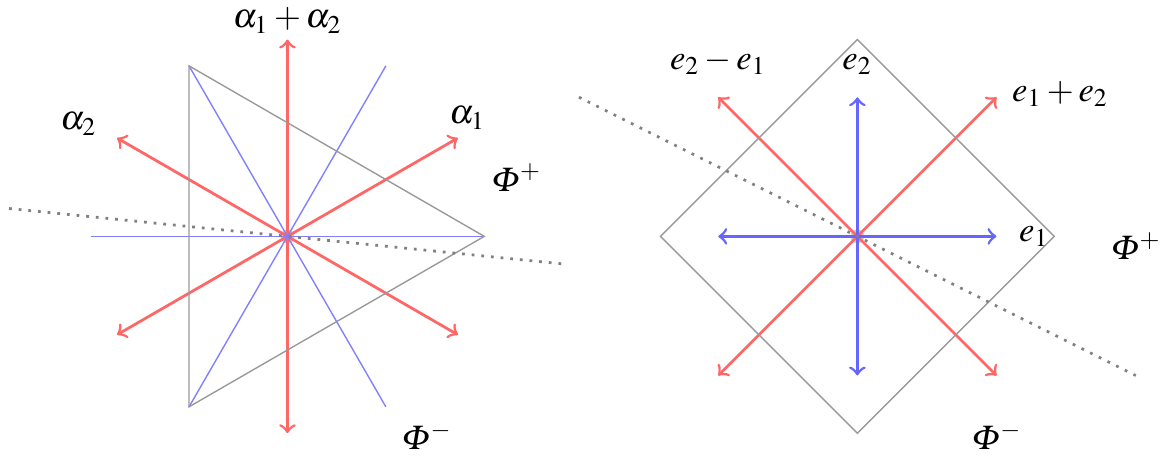}
   \caption{The picture on the  left  is a root system for the dihedral group $\mathcal D_3$ and the picture on the right is  a root system for the dihedral group $\mathcal D_4$.}
 \label{H_fig:rootdihedral} 
 \end{center}
 \end{figure}   %%%%%%%%%
\end{example}

\index{roots!positive roots}
\index{roots!negative roots}
\index{roots!simple roots}
Thanks to root systems, we are now able to find a suitable set of generators for $W$. The root system $\Phi$ is finite, so we can pick a hyperplane of $V$ which does not intersect $\Phi$ (represented with a dashed line in Figure~\ref{H_fig:rootdihedral}). This hyperplane induces a partition of $\Phi$ into two sets of the same cardinality, namely the set of {\em positive roots} $\Phi^+$ and the set of {\em negative roots} $\Phi^-=-\Phi^+$. It is obvious that $W$ is generated by the reflections $s_\alpha$ with $\alpha\in\Phi^+$. The polyhedral cone generated by $\Phi^+$ has a unique basis $\Delta\subseteq \Phi^+$ called a {\em simple system} consisting of roots  called {\em simple roots}. So any root is either a nonnegative linear combination of simple roots or a nonpositive linear combination of simple roots. In Figure~\ref{H_fig:rootdihedral}, the simple roots for $\mathcal D_3$ are $\alpha_1$ and $\alpha_2$, while for $\mathcal D_4$ they are $e_1$ and $e_2-e_1$.  We can thus set 
\bez
S=\Hset{s_\alpha}{s_\alpha\in\Delta}.
\eez
The reflections in $S$ are called {\em simple reflections}. The identity of $W$ is denoted by $e$.

\begin{theorem}[{see \cite[Chapter 1]{H_Hu90}}] 
\label{H_th:Coxeter}
Let $W$ be a finite reflection group with root system $\Phi$, set of positive roots $\Phi^+$ and simple system $\Delta$. 
\begin{enumerate}
\item The simple roots are linearly independent.
\item Any reflection $s$ in $W$ is conjugate in $W$ to a simple reflection. Moreover,  for a reflection $s$ in $W$, there is a unique  $\beta\in\Phi^+$ such that $s=s_\beta$.
\item $(W,S)$ is a {\em finite Coxeter system}: $W$ is generated by $S$ and by the  relations $s^2=e$ ($s\in S$ is a reflection) and $(st)^{o(st)}=e$ where $o(st)$ is the order of the rotation $st$, $s,t\in S$.
\end{enumerate}
\end{theorem}

\begin{remark}
\label{H_re:classi}
 The fact that finite reflection groups are finite Coxeter groups is crucial for their classification. Indeed, finite Coxeter groups are nicely classified through their root systems, and so are finite reflection groups, see N.~Reading's text in this volume~\cite{H_Re11}. 
\end{remark}

\begin{example}
\label{H_ex:A3Coxeter}
 For the symmetric group $S_n$, consider the set
\bez
\Phi=\Hset{e_j-e_i}{1\leq i\not = j\leq n}.
\eez
A positive root system of $\Phi$ is 
\bez
\Phi^+=\Hset{e_j-e_i}{1\leq i<j \leq n}.
\eez
A simple system is $\Delta=\Hset{e_{i+1}-e_i}{1\leq i < n}$ and the corresponding simple generators are the simple transpositions $\tau_i=(i\ i+1)$, well-known to generate $S_n$. 
\end{example}

\begin{example}
\label{H_ex:BnCoxeter}
 For the hyperoctahedral group $W'_n$, consider the set
\bez
\Phi=\Hset{e_j-e_i}{1\leq i\not = j\leq n}\cup\Hset{\pm e_i}{1\leq i\leq n}.
\eez
A positive root system of $\Phi$ is 
\bez
\Phi^+=\Hset{e_j-e_i}{1\leq i<j \leq n}\cup\Hset{e_i}{1\leq i\leq n}.
\eez
A simple system is $\Delta=\Hset{e_1,e_{i+1}-e_i}{1\leq i < n}$ and the corresponding simple generators are the simple transpositions $s_0$ and $\tau_i=(i\ i+1)$. 
\end{example}

\begin{example}
\label{H_ex:H3Coxeter}
 For the isometry group of the icosahedron  $W(H_3)$,  we have the set of simple generators $S=\{s_1,s_2,s_3\}$ with  $s_i^2=e$ and the other relations come from the order of rotations: $(s_1s_2)^5=(s_1s_3)^2=(s_2s_3)^3=e$. There are $15$ positive roots and a simple system consists of the following vectors in $\mathbb R^3$:
 \bez
   \left(2,0,0 \right),\qquad  \left(\frac{-1-\sqrt 5}{2},\frac{-1+\sqrt 5}{2},-1 \right),\qquad  \left(0,0,2 \right).
 \eez
\end{example}

\subsubsection{Permutahedra as $\mathcal H$-polytope} %%%%%%%
\label{H_ss:PermH}

We are now ready to complete our description of permutahedra as an $\mathcal H$-polytope. Fix a generic point $\pmb a$ and consider the associated permutahedron $\textnormal{Perm}^{\pmb a}(W)$. As we have seen in Example~\ref{H_ex:SnPerm}, a permutahedron can live in an affine subspace not containing the origin $0$. Let us first explain this phenomenon. The subspace $V_0$ of $V$ spanned by any root system $\Phi$ of $W$  is stable under the action of $W$, and so is its orthogonal complement\footnote{The orthogonal complement of $V_0$ is fixed pointwise by $W$ since it is the intersection of all the hyperplanes associated to reflections in $W$.}.   So the affine subspace $V_{\pmb a}= \pmb a + V_0$  of $V$ directed by $V_0$ and passing through the point $\pmb a$ is also stable under the action of $W$, and therefore
\bez
 \textnormal{Perm}^{\pmb a}(W) \subseteq V_{\pmb a}.
\eez

Comparing Figure~\ref{H_fig:permdihedral} and Figure~\ref{H_fig:rootdihedral}, we observe that the edges of $\textnormal{Perm}^{\pmb a}(W)$ are directed by roots. This phenomenon holds for permutahedra of higher dimension.  Therefore, the facets have to be directed by affine subspaces spanned by subsets of roots, and these subspaces have to be the boundaries of some halfspaces whose intersection is $\textnormal{Perm}^{\pmb a}(W)$. More precisely, 
choose the unique simple system\footnote{In fact, it is equivalent to choose a fundamental chamber in the Coxeter complex, which is always possible,  see \cite[\S1.12]{H_Hu90} for more details.} $\Delta$ of $\Phi$ such that $\Hscalprod{\pmb a}{\alpha}>0$ for all $\alpha\in \Delta$.  For $\alpha\in\Delta$ a simple root, we consider the halfspace in $V_0$
\bez 
\mathcal H_{0}(\alpha):=\Hset{\pmb x=\sum_{\beta\in\Delta} x_\beta \beta \in V_0}{x_\alpha \leq 0}
\eez
with its boundary 
$H_{0}(\alpha):=\textnormal{span}\,(\Delta\setminus\{\alpha \})\subseteq V_0$.
Then we consider their affine counterparts passing through the point $\pmb a$:
\bez
\mathcal H_{\pmb a}(\alpha)= \pmb a + \mathcal H_{0}(\alpha)\quad\textrm{and}\quad
    H_{\pmb a}(\alpha)= \pmb a + H_{0}(\alpha).
\eez

\begin{theorem}
\label{H_thm:PermH-rep}
 Let $\pmb a$ be a generic point in $V$.
 \begin{enumerate}
  \item The permutahedron $\textnormal{Perm}^{\pmb a}(W)$, as an $\mathcal H$-polytope, is given by
   \bez
        \textnormal{Perm}^{\pmb a}(W) = \bigcap_{{\alpha\in \Delta\atop w\in W}}w( \mathcal H_{\pmb a}(\alpha))\subseteq V_{\pmb a}.
   \eez
  \item For $w\in W$, $w(\pmb a)$ is a vertex of $\textnormal{Perm}^{\pmb a}(W)$ and 
  \bez
        \{w(a)\} = \bigcap_{{\alpha\in \Delta}}w(H_{\pmb a}(\alpha)).
   \eez
   \item The permutahedron $\textnormal{Perm}^{\pmb a}(W)$ is a $|\Delta|$-dimensional simple polytope. In particular, $\textnormal{Perm}^{\pmb a}(W)$ is full dimensional in the affine space $V_{\pmb a}$.
      \end{enumerate}
\end{theorem}

\begin{example} We illustrate this theorem in Figure~\ref{H_fig:PermDi4-roots}  with the example of the permutahedron of the dihedral group $\mathcal D_4$. 

 \begin{figure}  %%%%%%%%%Figure
 \begin{center}
  \includegraphics[scale=0.8]{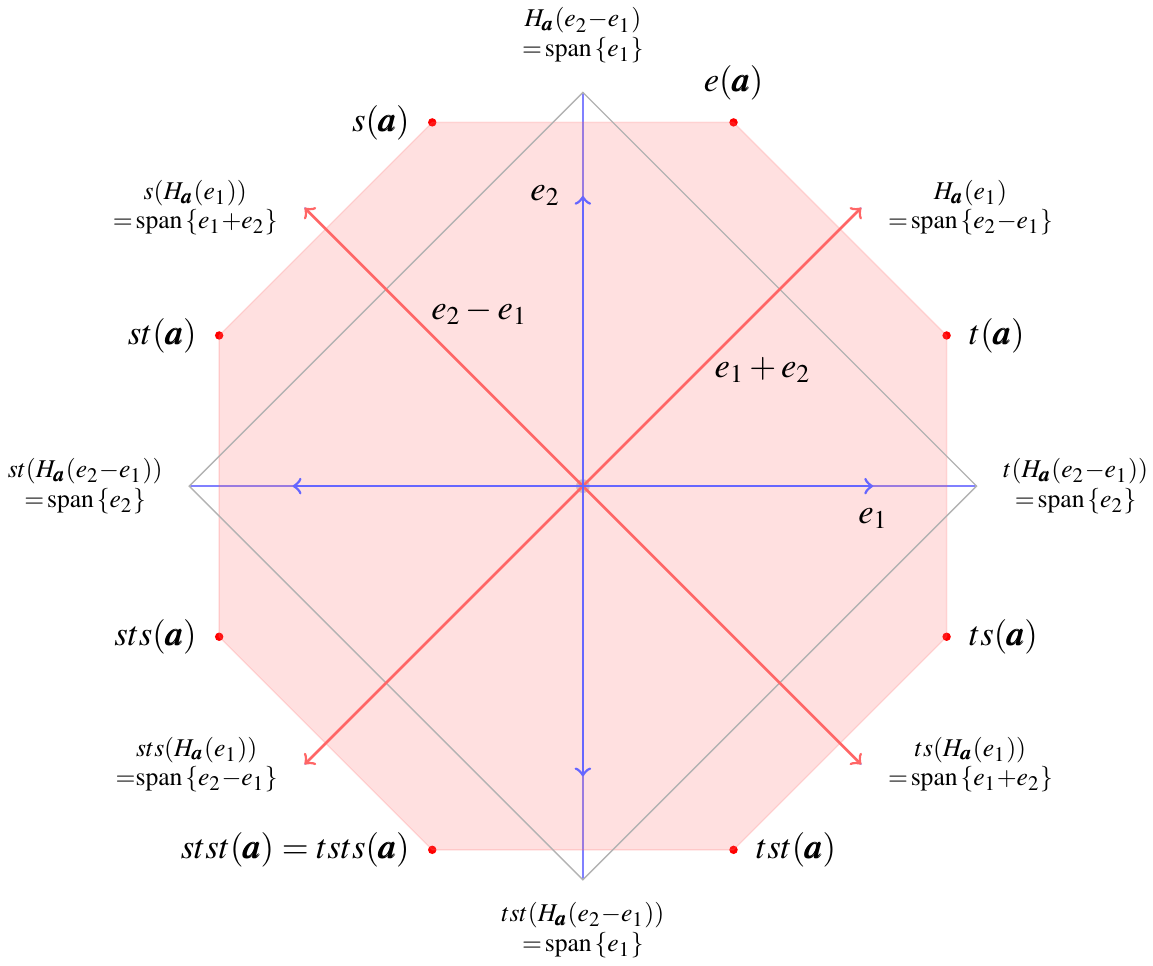}
   \caption{A  permutahedron  of the dihedral group $W=\mathcal D_4$ given as an $\mathcal H$-polytope. In this picture the roots are colored in blue and red as in Figure~\ref{H_fig:rootdihedral}; the simple reflections are $s:=s_{e_1}$ and $t:=s_{e_2-e_1}$, and  $V_0=V=V_{\pmb a}$.}
 \label{H_fig:PermDi4-roots} 
 \end{center}
 \end{figure}   %%%%%%%%%
\end{example}

The normal fan of a polytope $P$ is a collection of pointed polyhedral cones indexed by the faces of $P$: if $F$ is a face of $P$, the cone $C_F$ in the normal fan is the polyhedral cone generated by the outer normal vectors of the facets of $P$ containing $F$. Normal fans are natural objects that link polytope theory to optimization problems or algebraic geometry.  For more information on normal fans of polytopes, see \cite[Chapter 7]{H_Zi95}.

\begin{definition}
\index{Coxeter fan}
The normal fan $\mathcal F_W$ of $\textnormal{Perm}^{\pmb a}(W)$ is called the {\em Coxeter fan}.
\end{definition}

\begin{remark}[Note on the proof of Theorem~\ref{H_thm:PermH-rep}]  
As far as our knowledge goes, we do not know that many references which present the construction of permutahedra of finite reflection groups. Nevertheless, it was not a new idea. The proof is based on well-known and very important properties of finite reflection groups (see \cite[\S1.12]{H_Hu90}): the complement of the union of all hyperplanes corresponding to reflections in $W$ cut the space $V$ into open convex polyhedral cones called {\em chambers}. The collection of cones obtained by decomposing the boundaries of the chambers is the {\em Coxeter fan}. Picking one of these chambers to be the {\em fundamental chamber}, one can show that any chamber is the image of the fundamental chamber by an element of $W$ and that the isotropy group of a chamber is the identity. This implies that there are $|W|$ chambers and that any chamber contains exactly one point of the $W$-orbit of $\pmb a$.  Then observe that the hyperplanes $H_{\pmb a}(\alpha)$ are  the halfspaces orthogonal to the rays of the fundamental chambers passing through  $\pmb a$ and are stable under the subgroup $W_{S\setminus\{s_\alpha\}}$ generated by $\Hset{s_\beta}{\beta\in\Delta\setminus\{\alpha\}}$. Finally one shows, with a bit of straightforward work, that for any $\alpha \in \Delta$ and $w\in W\setminus W_{S\setminus\{s_\alpha\}}$,  the point $w(\pmb a)$ is in the interior of the halfspace $\mathcal H_{\pmb a}(\alpha)$.  For more details, see \cite{H_BoBo10,H_HoLaTh11}. 
\end{remark}

\subsection{Faces of permutahedra and the weak order} %%%%%%%%%%%%

Let $\pmb a$ be a generic point, $\Delta$ the simple system of $\Phi$ such that $\Hscalprod{\pmb a}{\alpha}>0$ for all $\alpha\in \Delta$ and $S=\Hset{s_\alpha}{\alpha\in\Delta}$ the set of simple reflections generating $W$. Our aim in this final part of our study of permutahedra is to explain how to realize several important notions from the theory of finite reflection groups in the context of permutahedra.

\subsubsection{Faces and standard parabolic subgroups}

The subgroup $W_I$ of $W$ generated by $I\subseteq S$ is called a {\em standard parabolic subgroup}. Since $W_I$ is generated by reflections, it is also a finite reflection group. Moreover, one can show that $\Delta_{I}:=\Hset{\alpha\in\Delta}{s_\alpha\in I}$ is a simple system for the root system $\Phi\cap\textnormal{span}\,(\Delta_I)$. So we can apply our construction of permutahedra to $W_I$: the polytope
\bez
 F_I:=\textnormal{Perm}^{\pmb a}(W_I)=\textnormal{conv}\,\Hset{w(\pmb a)}{w\in W_I}
\eez
is a $W_I$-permutahedron.  For instance, if $I=S\setminus\{s_\alpha\}$, then the hyperplane $H_{\pmb a}(\alpha)$ is stable under the action of $W_I$ and contains $F_I$, which is a facet of $\textnormal{Perm}^{\pmb a}(W)$. More generally, $F_I$ is a $|I|$-dimensional face of $\textnormal{Perm}^{\pmb a}(W)$ containing $\pmb a$ and lives in the affine space
\bez
\bigcap_{\alpha\in\Delta\setminus\Delta_I} H_{\pmb a}(\alpha) = \pmb a + \textnormal{span}\,(\Delta_I).
\eez
The next statement explains that the faces of $\textnormal{Perm}^{\pmb a}(W)$ are obtained as the $W$-orbit of the faces containing $\pmb a$.

\begin{proposition} 
\begin{enumerate}
\item Each face of dimension $k$ of $\textnormal{Perm}^{\pmb a}(W)$ is a permutahedron: $w(F_I)=\textnormal{Perm}^{w(\pmb a)}(W_I)$ for $w\in W$ and $I\subseteq S$ of cardinality $k$.  
\item For $w, g\in W$ and $I\subseteq S$, $w(F_I)=g(F_I)$ if and only if $wW_I=gW_I$.  In other words, faces are  naturally parametrized by the cosets $W/W_I$, $I\subseteq S$.
\item The face lattice of $\textnormal{Perm}^{\pmb a}(W)$ is isomorphic to the poset of cosets of standard parabolic subgroups: $w(F_I)\subseteq g(F_J)$ if and only if $wW_I\subseteq gW_J$.
\end{enumerate}
\end{proposition}

\begin{example}
We continue the example of the permutahedron of the dihedral group $\mathcal D_4$ we started in \S\ref{H_ss:perm}. Let $s=s_{e_1}$ and $t=s_{e_2-e_1}$; we are now able to label the vertices of $\textnormal{Perm}^{\pmb a}(\mathcal D_4)$ shown in Figure~\ref{H_fig:permdihedral}. For $w\in \mathcal D_4$ we label the vertex $w(\pmb a)$ by $w$. In particular the vertex $\pmb a$ is labeled by the identity $e$. As we can see in Figure~\ref{H_fig:PermDi4}, the facets, which are the edges here, are naturally labeled by the cosets of the subgroup $W_{s}$ generated by $s$ and the cosets of the subgroup $W_{t}$ generated by $t$. 

 \begin{figure}  %%%%%%%%%Figure
 \begin{center}
  \includegraphics[scale=0.8]{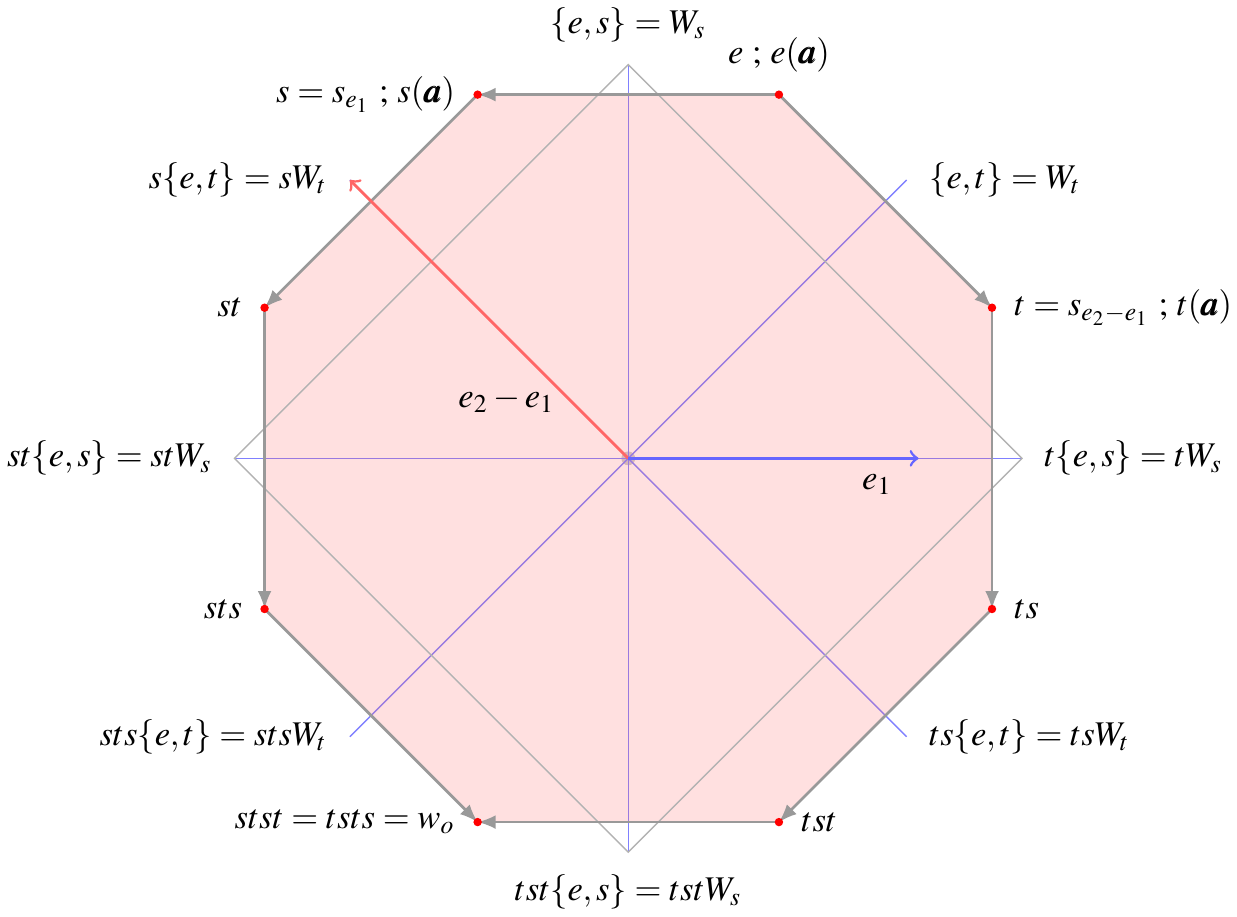}
   \caption{The labeled permutahedron  of the dihedral group $W=\mathcal D_4$ with nonempty faces labeled by the cosets of standard parabolic subgroups. The gray arrows, which are the oriented edges of the octagon, represent cover relations in the weak order of $W$.}
 \label{H_fig:PermDi4} 
 \end{center}
 \end{figure}   %%%%%%%%%
\end{example}

\subsubsection{Edges and weak order}

\begin{figure}  %%%%%%%%%Figure
 \begin{center}
  \includegraphics[scale=0.7]{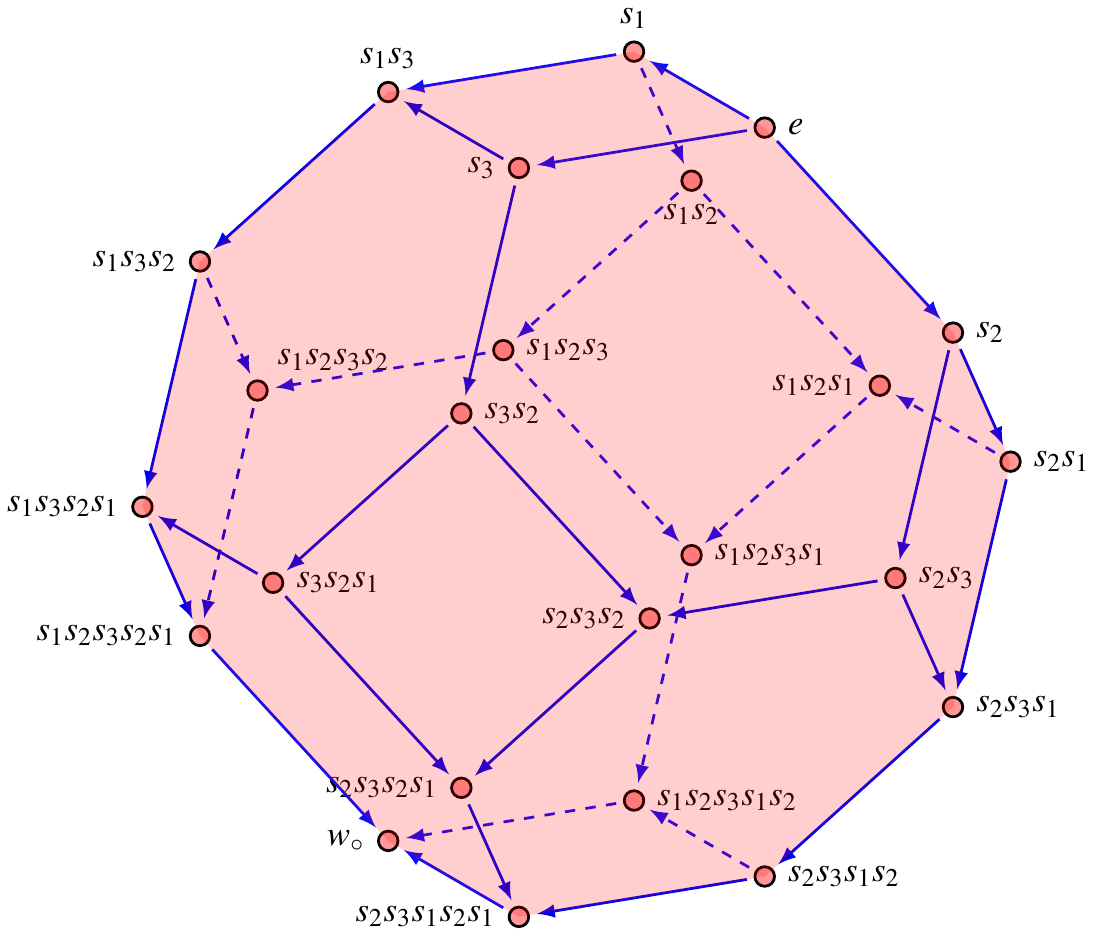}
   \caption{A permutahedron  of the group $S_4$ with labeled vertices and oriented edges which represent the weak order on $S_4$.}
 \label{H_fig:PermA3Weak} 
 \end{center}
 \end{figure}

We consider the length function $\ell:W\to\mathbb N$ mapping an element $w\in W$ to the minimal number $\ell(w)$ of letters needed to express $w$ as a word in the alphabet $S$. For instance $\ell(e)=0$ and $\ell(s)=1$ if and only if $s\in S$.  We know that an edge $w(F_{s})$ of $\textnormal{Perm}^{\pmb a}(W)$ corresponds to a coset $wW_{s}=\{w,ws\}$ for $s\in S$ and $w\in W$.  We orient each edge $w(F_s)$ from  $w$ to $ws$ if $\ell(w)<\ell(ws)$, and from $ws$ to $w$ otherwise. See Figure~\ref{H_fig:PermDi4} and Figure~\ref{H_fig:PermA3Weak} for examples.

    %%%%%%%%%

\begin{proposition}
The oriented $1$-skeleton\footnote{The $1$-skeleton of a polytope is the graph obtained by taking the vertices and the edges of this polytope.} of $\textnormal{Perm}^{\pmb a}(W)$ is the graph of the lattice called {\em the right weak order} of $W$. Moreover, the minimal element is  the identity $e$, the maximal element is denoted by $w_o$ and each face of $\textnormal{Perm}^{\pmb a}(W)$  corresponds to an interval in the weak order.
\end{proposition}

\begin{remark} 
\begin{enumerate}
\item The length $\ell(w_o)$ of $w_o$ is $|\Phi^+|$.  
\item The $f$-vector\footnote{The $k$-th coordinate  $f_k$ of the $f$-vector of a polytope $P$  is the number of $k$-dimensional faces of $P$, see \cite{H_Zi95}.}  of a permutahedron of $W$ encloses important information about the group: $f_0$ is the cardinality of $W$, $f_1$ is the number of cover relations in the weak order or $f_{n-1}$ is the number of parabolic subgroups of $W$.
\end{enumerate}
\end{remark}

\section{Coxeter generalized associahedra}

We are now ready to build Coxeter generalized associahedra. The basic idea is to build it as an $\mathcal H$-polytope by identifying  subsets of the halfspaces used in constructing a given permutahedron. These subsets of halfspaces will be chosen according to certain elements in $W$, namely {\em Coxeter singletons}, that are contained in the boundaries of the halfspaces. 

Throughout this section, we fix a generic point $\pmb a$, a simple system $\Delta$ such that $\Hscalprod{\pmb a}{\alpha}>0$ for all $\alpha\in \Delta$ and $S=\Hset{s_\alpha}{\alpha\in\Delta}$ the set of simple reflections which generates $W$. 

\subsection{Coxeter elements and Coxeter singletons}

Geometrically, Coxeter singletons correspond to vertices in a permutahedron of $W$  along certain paths between the identity $e$ and the longest element $w_o$ (a kind of spinal cord of the permutahedron). Let us be more precise. 

Let~$c$ be a Coxeter element of~$W$, that is, the product of the simple reflections in $S$ 
taken in some order, and fix a reduced expression for $c$ as a word in the alphabet $S$. For $I\subseteq S$, we denote by $c_{(I)}$ the subword of $c$ obtained by taking
only the simple reflections in $I$. So $c_{(I)}$ is a Coxeter element of $W_I$. For instance, if $W=S_4$ with $S=\{\tau_1,\tau_2,\tau_3\}$, a particular Coxeter element is $c=\tau_2\tau_3\tau_1$;  if  $I=\{\tau_1,\tau_2\}\subseteq S$, then $c_{(I)}=\tau_2\tau_1$.

 The longest element $w_o$ of the $W$ can be written as a reduced word on $S$ in many ways: each word corresponds precisely to a minimal path from $e$ to $w_o$ on the $1$-skeleton of a permutahedron for $W$. If we choose a Coxeter element $c$, N.~Reading showed in his work on {\em Coxeter sortable elements} that we can {\em sort} a particular reduced expression for $w_o$ accordingly to $c$ \cite{H_Re07-1}. This particular word is called {\em the $c$-word of $w_o$} and it is denoted by $\pmb{w_o}(c)$.  More precisely,  $\pmb{w_o}(c)$ is the unique reduced expression for $w_o$ on the alphabet\footnote{We have to make a clear distinction between words on the alphabet $S$ and reduced expressions of elements of $W$. The latter are subject to relations of the Coxeter system $(W,S)$ while the first ones are not.} $S$  such that  $\pmb{w_o}(c)=c_{(K_1)}c_{(K_2)}\dots c_{(K_p)}$ with non-empty~$K_i \subseteq S$, $K_p\subseteq K_{p-1}\subseteq\dots\subseteq K_1$ and  $\ell(w_o)=\sum_{i=1}^p |K_i|$. For instance, if $W=S_4$ with $S$ the set of simple transpositions $\tau_i$, we have the following $c$-words for $w_o$:
 \bez 
 \pmb{w_o}(\tau_1\tau_2\tau_3)= \tau_1\tau_2\tau_3.\tau_1\tau_2.\tau_1=c_{(S)}c_{(\{\tau_1,\tau_2\})}c_{(\{\tau_1\})}
 \eez
 and
 \bez
 \pmb{w_o}(\tau_2\tau_3\tau_1)= \tau_2\tau_3\tau_1.\tau_2\tau_3\tau_1 = c_{(S)}c_{(S)}.
 \eez
 
 \index{Coxeter singletons}
 \begin{definition} Let $c$ be a Coxeter element and $u\in W$. We say that  $u$ is a {\em $c$-singleton} if some reduced word $\pmb u(c)$ for $u$ appears as a prefix of a word that can be obtained from $\pmb{w_o}(c)$ by the commutation of commuting reflections of $S$.
The word $\pmb u (c)$ is called the {\em $c$-word} of the $c$-singleton $u$. 
 \end{definition}

\begin{example} Consider again the symmetric  group~$W=S_4$ together with $S=\{\tau_1,\tau_2,\tau_3\}$ the set of simple transpositions. The reader may follow this example with Figure~\ref{H_fig:AssA3-1}  and Figure~\ref{H_fig:AssA3-2} in mind for an illustration on a permutahedron.

For the Coxeter element $c=\tau_1\tau_2\tau_3$, the $c$-singletons, and their $c$-words, are
\begin{center}
\begin{tabular}{llll}
   $e$,            & $\tau_1$,             & $\tau_1\tau_2$, & $\tau_1\tau_2\tau_3$, \\
   $\tau_1\tau_2\tau_1,\qquad$ &  $\tau_1\tau_2\tau_3\tau_1,\qquad$ & $\tau_1\tau_2\tau_3\tau_1\tau_2,\qquad$ & $w_o=\tau_1\tau_2\tau_3\tau_1\tau_2\tau_1$.
\end{tabular}
\end{center} 
Observe that $\tau_1\tau_2\tau_1$ is a not a prefix of the word $\pmb{w_o}(c)=\tau_1\tau_2\tau_3\tau_1\tau_2\tau_1$, but it is a prefix of the word $\tau_1\tau_2\tau_1\tau_3\tau_2\tau_1$  obtained
after exchanging the commuting simple transpositions $\tau_1$, $\tau_3$ in the word $\pmb{w_o}(c)$.

For the Coxeter element $c'=\tau_2\tau_1\tau_3$, the $c'$-singletons, and their $c'$-words, are
\begin{center}
\begin{tabular}{lll}
   $e$,            & $\tau_2\tau_3$,             & $\tau_2\tau_3\tau_1\tau_2\tau_3$,\\
   $\tau_2$,          & $\tau_2\tau_3\tau_1$,          & $\tau_2\tau_3\tau_1\tau_2\tau_1$, and\\
   $\tau_2\tau_3,\qquad$ & $\tau_2\tau_3\tau_1\tau_2,\qquad$ & $w_o=\tau_2\tau_1\tau_3\tau_2\tau_1\tau_3$.
\end{tabular}
\end{center}
Observe that $\tau_2\tau_3$ is a not a prefix of the word $\pmb{w_o}(c')=\tau_2\tau_1\tau_3\tau_2\tau_1\tau_3$, but it is a prefix of the word $\tau_2\tau_3\tau_1\tau_2\tau_3\tau_1$  obtained after exchanging the
 commuting simple transpositions $\tau_1$, $\tau_3$ in the word $\pmb{w_o}(c')$.

\end{example}

\begin{problem} Let $W$ be a finite reflection group with set of simple reflections $S$. Find a formula for the number of $c$-singletons as a function of $c$. For which $c$'s are the maximum and the minimum reached? (As we can see in the examples above, the number of $c$-singletons depends on the choice of the Coxeter element $c$ or equivalently, on the chosen orientation of the Coxeter graph, see Remark~\ref{H_re:orientation} below. Numerical examples are presented in~\cite{H_HoLa07}). 
\end{problem}

\begin{remark}
\label{H_re:orientation}
\index{Coxeter sortable elements}
\index{Tamari lattice}
\index{Cambrian lattices}
\begin{enumerate}
\item  {\em Coxeter sortable elements.} The set of elements in $W$ that verify the same characterization we gave above for the $c$-word $\pmb{w_o}(c)$ of the longest element $w_o$ are called {\em $c$-sortable elements}. They were introduced by N.~Reading for studying extensions of the {\em Tamari lattice} called {\em Cambrian lattices}~\cite{H_Re07-1,H_Re07-2}. The combinatorics of these elements is very rich (see for more details N.~Reading's text in this volume~\cite{H_Re11}).
\item Coxeter elements are in bijection with {\em orientations of Coxeter graphs}. A Coxeter graph is a graph whose vertices are simple reflections and whose edges are labeled by the order of the product of two simple reflections (a rotation); there is no edge between two commuting simple reflections.  Any Coxeter element~$c$ defines an orientation of the Coxeter graph of $W$:
orient the edge $\{s_i,s_j\}$ from~$s_i$ to~$s_j$ if and only if~$s_i$ is to the left of~$s_j$ in any reduced word for~$c$. The contents of the articles \cite{H_HoLa07,H_BeHoLaTh11,H_HoLoRa10} makes use of this bijection, as well as N.~Reading in his works on Coxeter sortable elements (see~\cite{H_Re11}).  
\end{enumerate}
\end{remark}

\subsection{Coxeter generalized associahedra as $\mathcal H$-polytopes}

At this point, one should remember that the permutahedron obtained from the generic point $\pmb a$ is described as an $\mathcal H$ polytope by
\bez
              \textnormal{Perm}^{\pmb a}(W)=\bigcap_{\alpha\in \Delta,\atop\ s \in W} w\big(\mathcal H_{\pmb a}(\alpha)\big)\subseteq V_{\pmb a}\quad \textrm{(see \S\ref{H_ss:PermH})}.
            \eez

\begin{definition}
\index{associahedron!Coxeter generalized associahedra}
\label{H_def:Asso}
 Let~$c$ be a Coxeter element of~$W$. The polytope in the affine space $V_{\pmb a}$ defined by
            \bez
              \textnormal{Asso}^{\pmb a}_c(W)=\bigcap_{\alpha\in \Delta,\atop\ u \textrm{ is a $c$-singleton}} u\big(\mathcal H_{\pmb a}(\alpha)\big)
            \eez
 is called a {\em $c$-generalized associahedron}. As we will see, the combinatorics of $\textnormal{Asso}^{\pmb a}(W)$ does not depend of the choice of $\pmb a$, as long as this point is generic.
\end{definition}

As a first consequence, we see that for $\alpha\in\Delta$ and $u\in W$ a $c$-singleton, the halfspace $u(\mathcal H_{\pmb a}(\alpha))$ contains the point $u(\pmb a)$ in its boundary hyperplane $u(H_{\pmb a}(\alpha))$. We say that the halfspace $w(\mathcal H_{\pmb a}(\alpha))$ is  {\em $c$-admissible} if its boundary $w(H_{\pmb a}(\alpha))$ contains a $c$-singleton. So we can see a Coxeter generalized associahedron as obtained from a given permutahedron by keeping only  the $c$-admissible halfspaces.

\begin{example}
\label{H_ex:AssD4}
\index{associahedron!Coxeter generalized associahedra! of dimension $2$}
 Coxeter generalized associahedra of dimension $2$ arise from dihedral groups: an associahedron for $\mathcal D_n$ is a $n+2$-gon. 

For instance, take $W=\mathcal D_4$, the symmetry group of a square as in Example~\ref{H_ex:PermD4}: take $s:=s_{e_1}$ and $t:=s_{e_2-e_1}$. Here $V_0=V=V_{\pmb a}$ . Consider the Coxeter element $c=ts$, then the $c$-singletons are $e$, $t$, $ts$, $tst$, $tsts$.   Figure~\ref{H_fig:AssoDi4} shows the resulting associahedron $\textnormal{Asso}^{\pmb a}(\mathcal D_4)$ containing the permutahedron $\textnormal{Perm}^{\pmb a}(\mathcal D_4)$.  Observe that the common vertices between $\textnormal{Perm}^{\pmb a}(\mathcal D_4)$ and  $\textnormal{Asso}^{\pmb a}(\mathcal D_4)$ are precisely the $c$-singletons. 
\begin{figure}  %%%%%%%%%Figure
 \begin{center}
  \includegraphics[scale=0.8]{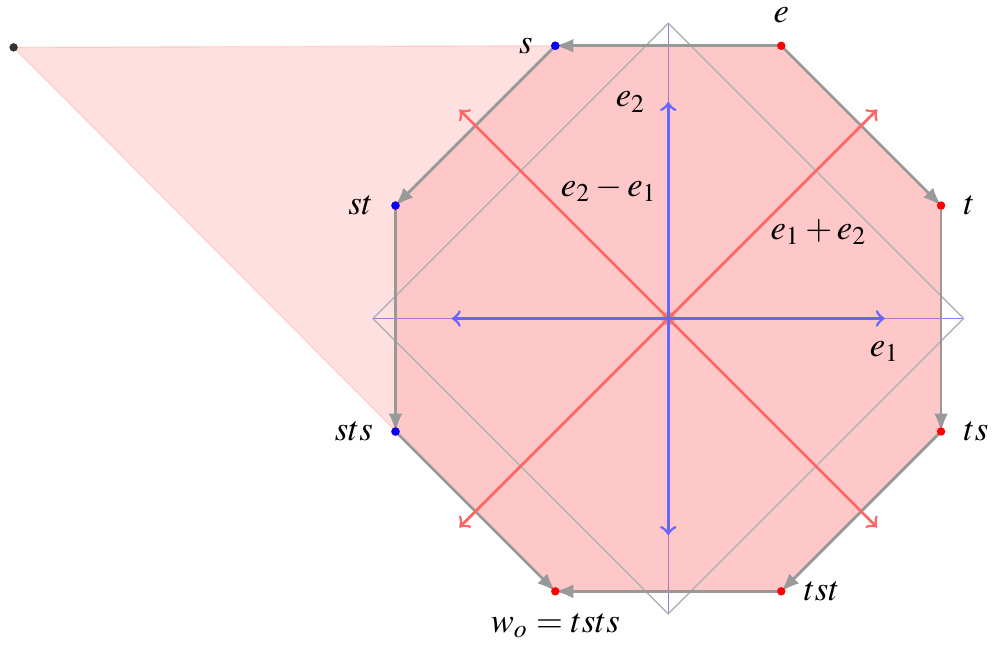}
   \caption{A Coxeter generalized  associahedron  of the dihedral group $W=\mathcal D_4$ given as an $\mathcal H$-polytope and containing the permutahedron with edges oriented by the weak order. In this picture, the Coxeter singletons are in red and the boundary of the $c$-admissible halfspaces meet in a vertex colored black}
 \label{H_fig:AssoDi4} 
 \end{center}
 \end{figure}   %%%%%%%%%
\end{example}

 \index{associahedron!Coxeter generalized associahedra! of dimension $3$}
\begin{example} There are two kinds of Coxeter generalized associahedra of dimension $3$.  
\label{H_ex:AssDim3}
The first kind is obtained in $\mathbb R^3$  from permutahedra of dimension $2$ by considering the isometry group of a regular polygonal prism as in Example~\ref{H_ex:PermDim3}. Consider a regular $n$-gonal prism with isometry group  $W=\mathcal D_n \times \mathcal D_2$. The Coxeter singletons are the couples $(u,e)$ and $(u,s)$ where $u$ is a Coxeter singleton of $\mathcal D_n$.  Then Coxeter generalized associahedra obtained  from $\textnormal{Perm}^{\pmb a}(W)$ are $n+2$-gonal prisms.

 \index{associahedron!classical associahedron}
\index{associahedron!cyclohedron}
The second kind of Coxeter generalized associahedra of dimension $3$ arises from the symmetric group $S_4$ (one of those is the {\em classical associahedron} also called {\em Stasheff polytope}), the hyperoctahedral group $W'_3$ (those are called {\em cyclohedron}) and the isometry group of the dodecahedron $W(H_3)$.  Examples are shown in  Figures~\ref{H_fig:AssA3-1},~\ref{H_fig:AssA3-2},~\ref{H_fig:AssB3},~\ref{H_fig:AssH3-1}~and~\ref{H_fig:AssH3-2}. 
\end{example}
 
The observations we can make on these pictures are summarized in the following theorem.

\index{associahedron!Coxeter generalized associahedra! as $\mathcal H$-polytope}
\begin{theorem}[\cite{H_HoLaTh11}]\label{H_thm:AssoC}
   Let~$c$ be a Coxeter element of~$W$.
   \begin{enumerate}
      \item $\textnormal{Asso}^{\pmb a}_c(W)$ is a $|\Delta|$-dimensional simple convex polytope.
      \item $\textnormal{Perm}^{\pmb a}(W)\subseteq \textnormal{Asso}^{\pmb a}_c(W)$.
      \item Each facet of $\textnormal{Asso}^{\pmb a}_c(W)$ is contained in the boundary of exactly one $c$-admissible halfspace. There are $(|\Delta|+|\Phi^+|)$ facets.
      \item The vertex sets~$\textnormal{vert}\,(\textnormal{Asso}^{\pmb a}_{c}(W))$ and~$\textnormal{vert}\,(\textnormal{Perm}^{\pmb a}(W))$ satisfy
            \bez
              \textnormal{vert}\,(\textnormal{Asso}^{\pmb a}_{c}(W))\cap \textnormal{vert}\,(\textnormal{Perm}^{\pmb a}(W))
                    =\{ u(\pmb a) \, | \, u \text{ is a $c$-singleton}\};
            \eez
            and this intersection forms a distributive sublattice of the weak order.
\end{enumerate}
\end{theorem}

\begin{definition}
\index{Cambrian fans}
The normal fan of $\textnormal{Asso}^{\pmb a}_c(W)$ is called the {\em $c$-Cambrian fan}.
\end{definition}

\begin{remark}[{Note on the proof of Theorem~\ref{H_thm:AssoC}}]
\index{associahedron!Coxeter generalized associahedra! as $\mathcal V$-polytope}
\index{spherical subword complex}
 The original motivation for constructing $\textnormal{Asso}^{\pmb a}_c(W)$ was to show that {\em Cambrian fans} are normal fans of some polytopes, answering a conjecture made by N.~Reading in~\cite{H_Re06}. Cambrian fans are defined as coarsening fans of the Coxeter fan, see \cite{H_ReSp09,H_Re11}. The proof the authors gave in \cite{H_HoLaTh11} is very technical and based on some technical properties of Cambrian fans given by N.~Reading and D.~Speyer. The principal reason for that difficulty is that we don't have a description of  $\textnormal{Asso}^{\pmb a}_c(W)$ as a $\mathcal V$-polytope, that is, as the convex hull of a given set of points. Recently, we were made aware of a very nice work in progress by V.~Pilaud and C.~Stump~\cite{H_PiSt11}: they introduce new families of polytopes called {\em brick polytopes} associated to finite reflection groups. A brick polytope is defined as a $\mathcal V$-polytope associated to a {\em spherical subword complex} of $W$, which were introduced by A.~Knutson and E.~Miller in~\cite{H_KnMi05}. They show in particular that brick polytopes contain Coxeter generalized associahedra providing a new proof of this theorem  and a description of $\textnormal{Asso}^{\pmb a}_c(W)$  as a $\mathcal V$-polytope.
\end{remark}

\subsection{Faces and almost positive roots}

The question now is to find a nice parameterization of the faces of $\textnormal{Asso}^{\pmb a}_c(W)$. 

\index{roots!almost positive roots}
Let $\Phi$ be a root system for $W$, with simple system $\Delta$. The set of almost positive roots is the set
\bez
  \Phi_{\geq -1} :=-\Delta\cup \Phi^+.
\eez

As stated in Theorem~\ref{H_thm:AssoC}, there are precisely $|\Phi_{\geq -1}|=|\Delta|+|\Phi^+|$ facets. This suggests a labeling of the facets by almost positive roots. To describe this labeling of the set of facets, we define for $\alpha\in\Delta$ the {\em last root map}  $\textnormal{lr}_\alpha$, which sends a $c$-singleton $u$ to an almost positive root $\textnormal{lr}_\alpha(u)\in\Phi_{\geq -1}$, as follows:
\begin{enumerate}
\item if $s_\alpha\in S$ is not a letter in the $c$-word $\pmb u(c)$ of $u$, then 
\bez 
  \textnormal{lr}_\alpha(u):=-\alpha\in - \Delta;
\eez
\item if $s_\alpha\in S$ is a letter in the word $\pmb u(c)$, we write $\pmb u(c) = u_1 s_\alpha u_2$ where $u_2$ is the unique largest suffix, possibly empty, of $\pmb u(c)$ that does not contain the letter $s_\alpha$ and 
\bez
  \textnormal{lr}_\alpha(u)=u_1 (\alpha)\in\Phi^+ .
\eez
\end{enumerate}

\begin{example}
      To illustrate this map, we consider again the dihedral group~$W=\mathcal D_4$ with
        generators $S=\{s,t\}$ as in Example~\ref{H_ex:AssD4}.
        Fix~$c=ts$ as a Coxeter element. We saw in Example~\ref{H_ex:AssD4}
        that the $c$-singletons are: $e$, $t$, $c=ts$ and $w_o=tsts$. 
       For the $c$-word $\pmb{w_o}(c)= tsts$ and $t=s_{e_2-e_1}$ we have 
       $\pmb{w_o}(c)=u_1t u_2$ with $u_1=ts$ and $u_2=s$. We get $\textnormal{lr}_{e_2-e_1}(tsts)=ts(e_2-e_1)=e_1+e_2$. 
       We show in Figure~\ref{H_fig:PermDi4-roots} the last root map $\textnormal{lr}_\alpha(u)$ for all the $c$-singletons $u$ and simple roots $\alpha$.
\end{example}

Thanks to the maps $\textnormal{lr}_\alpha$, we are able, together with Theorem~\ref{H_thm:AssoC}, to label the facets of $\textnormal{Asso}^{\pmb a}_c(W)$: for $\alpha\in\Phi_{\geq -1}$ denote by $F_{\textnormal{lr}_\alpha(u)}$ the unique facet of $\textnormal{Asso}^{\pmb a}_c(W)$ supported by the hyperplane $u(H_{\pmb a}(\alpha))$, where $\alpha\in \Delta$ and $u$ is a $c$-singleton . Therefore, since $\textnormal{Asso}^{\pmb a}_c(W)$ is a $|\Delta|$-dimensional simple polytope and since any face of codimension $k$ is the intersection of $k$ facets, we obtain a natural labeling of the faces of $\textnormal{Asso}^{\pmb a}_c(W)$: any face of $\textnormal{Asso}_c^{\pmb a}(W)$ of codimension $k$ is of the form $F_\Lambda$ where  $\Lambda=\{\alpha_1,\dots,\alpha_k\}\subseteq\Phi_{\geq -1}$ such that
\bez
  F_\Lambda=\bigcap_{i=1}^k F_{\alpha_{i}}.
\eez

\begin{example} This labeling is shown on Figure~\ref{H_fig:AssDi4-roots} for the example of the dihedral group~$W=\mathcal D_4$.
 \begin{figure}  %%%%%%%%%Figure
 \begin{center}
   \includegraphics[scale=0.85]{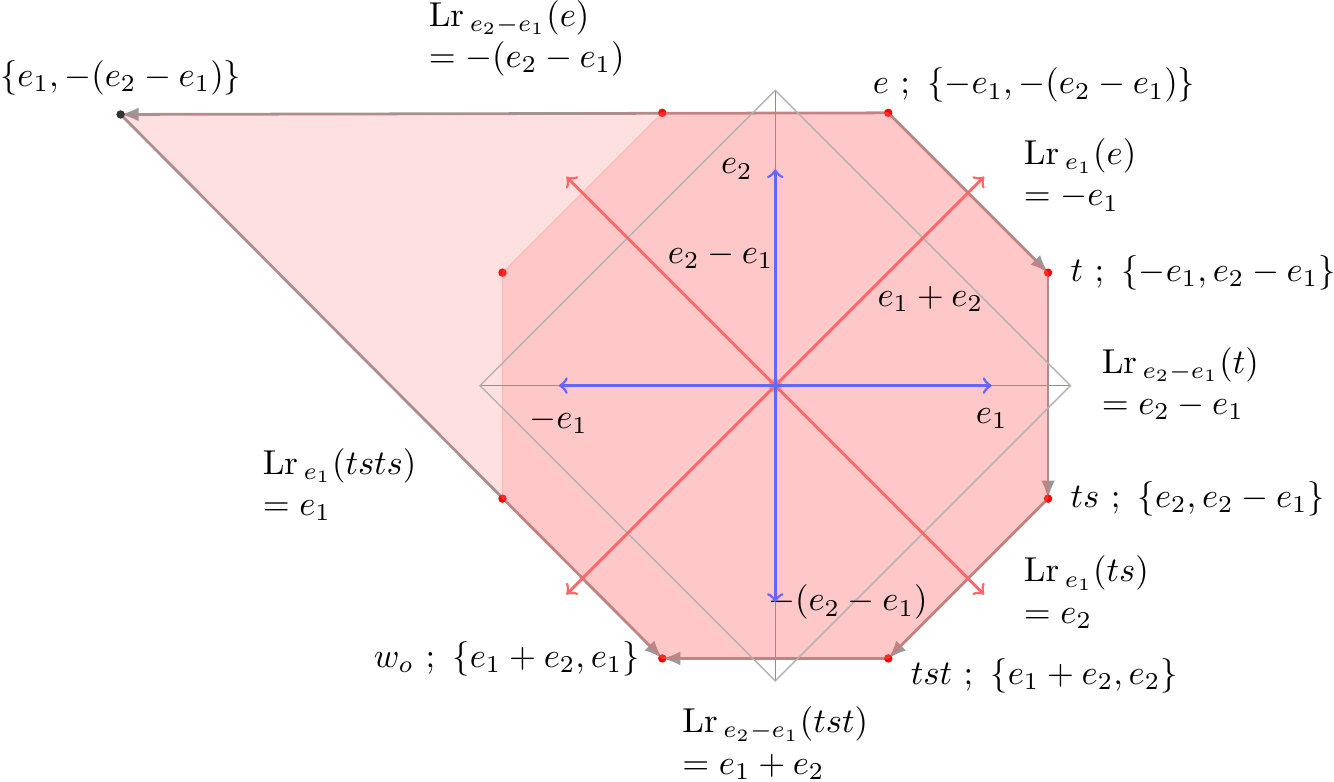}
   \caption{A generalized  associahedron  of the dihedral group $W=\mathcal D_4$ whose faces are labeled by subsets of almost positive roots. In this picture, the Coxeter singletons are in red. The edges are oriented according to the $c$-Cambrian lattice, see \S\ref{H_ss:Camblatt}.}
 \label{H_fig:AssDi4-roots} 
 \end{center}
 \end{figure}   %%%%%%%%%
\end{example}

\index{associahedron!Coxeter generalized associahedra!faces}
\begin{theorem} 
\label{H_thm:FaceAsso}
\begin{enumerate}
 \item The map $\beta\mapsto F_\beta$ is a bijection between the set of almost positive roots $\Phi_{\geq -1}$ and the facets of $\textnormal{Asso}^{\pmb a}_c(W)$.
\item The vertices of $\textnormal{Asso}^{\pmb a}_c(W)$ are  labeled subsets of $\Phi_{\geq -1}$ of cardinality $|\Delta|$, called {\em $c$-clusters} . The $1$-skeleton of $\textnormal{Asso}^{\pmb a}_c(W)$ is called  {\em the $c$-cluster exchange graph}. They are counted by the {\em $W$-Catalan numbers}.
\end{enumerate}
\end{theorem}

\index{clusters}
\index{Catalan numbers}
Catalan numbers appear in the context of symmetric groups. They count many objects such as planar binary trees, triangulations of a given polygon, noncrossing partitions, etc. They have an analog for any finite reflection group where they count, for instance,  clusters of finite type, Coxeter sortable elements or generalized noncrossing partitions, see for instance \cite{H_Re07-1} for more details. One way to prove that $c$-clusters are counted by the $W$-Catalan numbers is given by N.~Reading in \cite{H_Re07-1}: he constructs a bijection from $c$-clusters to noncrossing partitions via $c$-sortable elements.  For more on the connection to cluster algebras, see N. Reading's text in this volume~\cite{H_Re11}. 

\begin{problem} Find from the vertices of a $c$-generalized associahedron a direct and uniform proof  that the number of vertices is the $W$-Catalan number.
\end{problem}

\begin{remark}[Notes on the proof of Theorem~\ref{H_thm:FaceAsso}]
The first point of this theorem was proved in \cite{H_HoLaTh11}.  

In this same article, the authors showed the second point of the theorem by showing that the normal fan of this polytope is the {\em $c$-Cambrian fan} studied by N.~Reading and D.~Speyer in~\cite{H_ReSp09}.  N.~Reading and D.~Speyer showed in particular, using results obtained by N.~Reading in \cite{H_Re07-1,H_Re07-2}, that the maximal cones of the $c$-Cambrian fan are the $c$-clusters and that the graph obtained by considering maximal cones and walls in this fan is the $c$-cluster exchange graph. 

\index{spherical subword complex}
Recently, this connection between $c$-clusters and the vertices of $\textnormal{Asso}^{\pmb a}_c(W)$ was made entirely clear from a combinatorial point of view by C.~Ceballos, J.-P.~Labb\'e and C.~Stump in~\cite{H_CeLaSt11}. The authors find the $1$-skeleton of $\textnormal{Asso}^{\pmb a}_c(W)$ as the {\em facet-adjency graph of spherical subword complex} of the word $c\pmb{w_o}(c)$. This new object allow them, not only to recover the $1$-skeleton  of $\textnormal{Asso}^{\pmb a}_c(W)$ with the parameterizations above, but also to easily and naturally prove that the vertices are $c$-clusters in the sense of N.~Reading.
\end{remark}

\begin{remark}[The classical associahedron and the cyclohedron as Coxeter generalized associahedra]
\index{associahedron!classical associahedron}
\index{associahedron!cyclohedron}
\label{H_ex:HoLa}
In \cite{H_HoLa07}, the authors present Coxeter generalized associahedra for symmetric groups and hyperoctahedral groups both as $\mathcal V$-polytopes and $\mathcal H$-polytopes.  In that article, they start from the classical permutahedron $\Pi_{n-1}$.  For the symmetric group $S_n$, the faces of generalized associahedra can be labeled by triangulations of a $n+2$-gon $P$. The idea is to index the vertices of $P$ according to the chosen Coxeter element $c$, or equivalently to the orientation of the  Coxeter graph of $S_n$, see Remark~\ref{H_re:orientation}, and then to give a very easy combinatorial way to associate {\em integer coordinates} in $\mathbb R^n$ to any triangulation of $P$. These coordinates turn out to be precisely the coordinates of the vertices of the corresponding $c$-generalized associahedron of $S_n$.  This particular class of generalized associahedra, obtained from the classical permutahedron, contains the `classical' associahedron as constructed by J.-L.~Loday~\cite{H_Lo04} and S.~Shnider-S.~Sternberg~\cite{H_ShSt94}, and the associahedron arising from cluster algebra  theory as realized by F.~Chapoton, S.~Fomin and A.~Zelevinsky~\cite{H_ChFoZe03} (see Figures~\ref{H_fig:AssA3-1}~and~\ref{H_fig:AssA3-2}). In regard to the last statement, we refer the reader to the very nice article by C.~Ceballos, F.~Santos and G.~M.~Ziegler~\cite{H_CeSaZi11} on the subject of comparing these realizations.  Many realizations of the cyclohedron are easily obtained from this class of generalized associahedra of symmetric groups by looking at the hyperoctahedral group $W'_n$ as a subgroup of $S_{2n}$ (see~\cite{H_HoLa07} for more details).

Using this description as a $\mathcal V$-polytope, J.~Lortie, A.~Raymond and the author were able to show that the centers of the gravity of the vertices of the classical permutahedron and of all the Coxeter generalized associahedron built from it are the same~\cite{H_HoLoRa10}. This is still an open problem for an arbitrary Coxeter generalized associahedra.

\begin{problem} Let $\pmb a$ be a generic point such that $\Hscalprod{\pmb a}{\alpha}=\Hscalprod{\pmb a}{\beta}$ for all $\alpha,\beta\in \Delta$ and $W$ a finite reflection group. Prove that the vertices of any $c$-generalized associahedra $\textnormal{Asso}_c^{\pmb a}(W)$ and of  the permutahedron $\textnormal{Perm}^{\pmb a}(W)$  have the same center of gravity. (There is significant computational evidence supporting this statement).
\end{problem}

The proof in \cite{H_HoLoRa10} was based on the action of $\mathcal D_{n+2}$ on the set of vertices of the associahedron indexed by triangulations of $P$. It turns out that for each orbit, the center of gravity of the vertices in this orbit is also the same as that of the classical permutahedron. It leads us to ask the following question.

\begin{problem} For an arbitrary finite reflection group $W$, find an action on the vertices of a Coxeter generalized associahedron $\textnormal{Asso}_c^{\pmb a}(W)$ that generalizes the action of $\mathcal D_{n+2}$ on a regular $n+2$-gon, and such that the center of gravity of the vertices in each of the orbits is the same as  that of  $\textnormal{Perm}^{\pmb a}(W)$ if $\pmb a$ is a generic point verifying $\Hscalprod{\pmb a}{\alpha}=\Hscalprod{\pmb a}{\beta}$ for all $\alpha,\beta\in \Delta$.
\end{problem}
\end{remark}

\subsection{Edges and Cambrian lattices}
\label{H_ss:Camblatt}
\index{Cambrian lattices}
\index{Tamari lattice}

Cambrian lattices were introduced by N.~Reading in \cite{H_Re06}, see also N.~Reading's text in this volume~\cite{H_Re11}. The $c$-Cambrian lattice is a sublattice, as well as a quotient lattice, of the weak order of $W$. These lattices extend to the framework of finite reflection groups the celebrated {\em Tamari lattice} arising in the framework of symmetric groups. To be more precise:  the sets of {\em $c$-sortable elements}, introduced in Remark~\ref{H_re:orientation}, form a sublattice of the weak order called the {\em $c$-Cambrian lattice}.  We will now explain how to recover this order on the $1$-skeleton of $\textnormal{Asso}^{\pmb a}_c(W)$. The problem is to orient each edge $[F_{\lambda},F_{\mu}]$ of $\textnormal{Asso}_c^{\pmb a}(W)$, where $\lambda,\mu$ are $c$-clusters. Unfortunately, it is not as combinatorially easy as it was for recovering the weak order from the permutahedron. 

In the case where the  edge $[F_{\lambda},F_{\mu}]=[u(\pmb a),v(\pmb a)]$ with $u$, $v$ are $c$-singletons, we just keep the orientation given by the weak order. In the case where one of the vertices of $[F_{\lambda},F_{\mu}]$ is not a $c$-singleton we have to use the $c$-cluster map define by N.~Reading in \cite{H_Re07-1}. The idea is that, similar to $c$-singletons, any $w\in W$ has a $c$-word $\pmb{w}(c)$. We are thus able to extend the last root map $\textnormal{lr}$ to $W$ and then to define the {\em $c$-cluster map} $\textnormal{cl}_c$ from $W$ to the subsets of $\Phi_{\geq -1}$  as follows:
\bez
w\in W\mapsto \textnormal{cl}_c(w):=\Hset{\textnormal{lr}_\alpha(w)}{\alpha\in\Delta}\subseteq\Phi_{\geq -1}.
\eez
It turns out that N.~Reading shows that the restriction of this map to $c$-sortable elements is a bijection with the set of $c$-clusters. We are now able to orient the edges: orient the  edge $[F_{\lambda},F_{\mu}]$ from $F_{\lambda}$ to $F_{\mu}$ if $\textnormal{cl}_c^{-1}(\lambda)$ is smaller than $\textnormal{cl}_c^{-1}(\mu)$ in the weak order.
 We have to note that this way of orienting the edges of $\textnormal{Asso}^{\pmb a}_c(W)$ is not convenient. 
 
 \begin{problem} Find a combinatorial way to orient the edges of $\textnormal{Asso}^{\pmb a}_c(W)$ to recover the $c$-Cambrian lattice without the  use of the cluster map $\textnormal{cl}_c$.
\end{problem}

A stronger statement would be to answer the following problem.

\begin{problem} Find a combinatorial way to label the vertices of $\textnormal{Asso}^{\pmb a}_c(W)$ by $c$-sortable elements without the  use of the cluster map $\textnormal{cl}_c$.
\end{problem}

\begin{proposition}[\cite{H_Re07-2,H_HoLaTh11}] 
The oriented $1$-skeleton of $\textnormal{Asso}^{\pmb a}_c(W)$ is the graph of the $c$-Cambrian lattice.  Moreover, the $c$-singletons form a distributive sublattice of the $c$-Cambrian lattice: the minimal element is  the identity $e$, the maximal element is  $w_o$ and each face of $\textnormal{Asso}_c^{\pmb a}(W)$  corresponds to an interval in the $c$-Cambrian lattice.
\end{proposition}

Examples are shown in  Figures~\ref{H_fig:AssDi4-roots},~\ref{H_fig:AssA3-1}~and~\ref{H_fig:AssA3-2}.

\begin{remark} As unoriented graphs,  $1$-skeletons of different $c$-generalized associahedra $\textnormal{Asso}^{\pmb a}_c(W)$ are isomorphic, which implies that they have the same combinatorial type. Nevertheless, in general,  $c$-Cambrian lattices are not lattice isomorphic, and $c$-generalized associahedra are not isometric.
\end{remark}

\subsection{Isometry classes} %%%%%%%%%%%%%
\index{associahedron!isometry classes}

\begin{figure}
\begin{center}
  \includegraphics[scale=0.8]{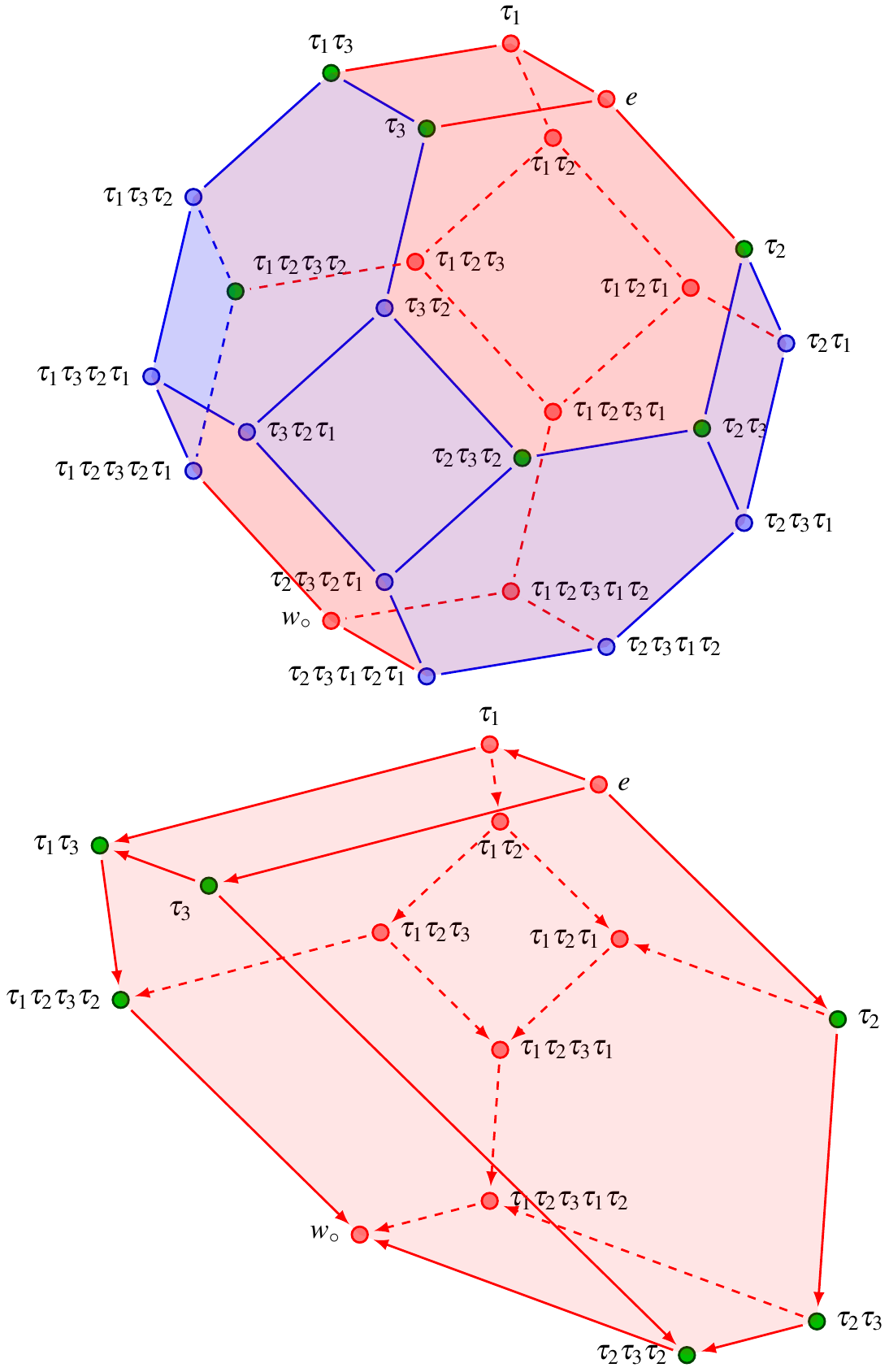}
\caption{{\bf The associahedron as realized by J.-L.~Loday and S.~Shnider-S.~Sternberg.} On the top is the classical permutahedron for the symmetric group $S_4$, see Example~\ref{H_ex:A3Coxeter},  and on the bottom is associahedron $\textnormal{Asso}^{\pmb a}_{\tau_1\tau_2\tau_3}(S_4)$ which corresponds to the associahedron realized by Loday and Shnider-Sternberg.  The facets of $\textnormal{Asso}^{\pmb a}_{\tau_1\tau_2\tau_3}(S_4)$ containing Coxeter singletons, as well as the  Coxeter singletons, are in red. In green are the vertices corresponding to the $\tau_1\tau_2\tau_3$-sortable elements. The orientations on the edges represent the $\tau_1\tau_2\tau_3$-Cambrian lattice which, is the Tamari lattice.}
\label{H_fig:AssA3-1}
\end{center}
\end{figure}

\begin{figure}
\begin{center}
  \includegraphics[scale=0.8]{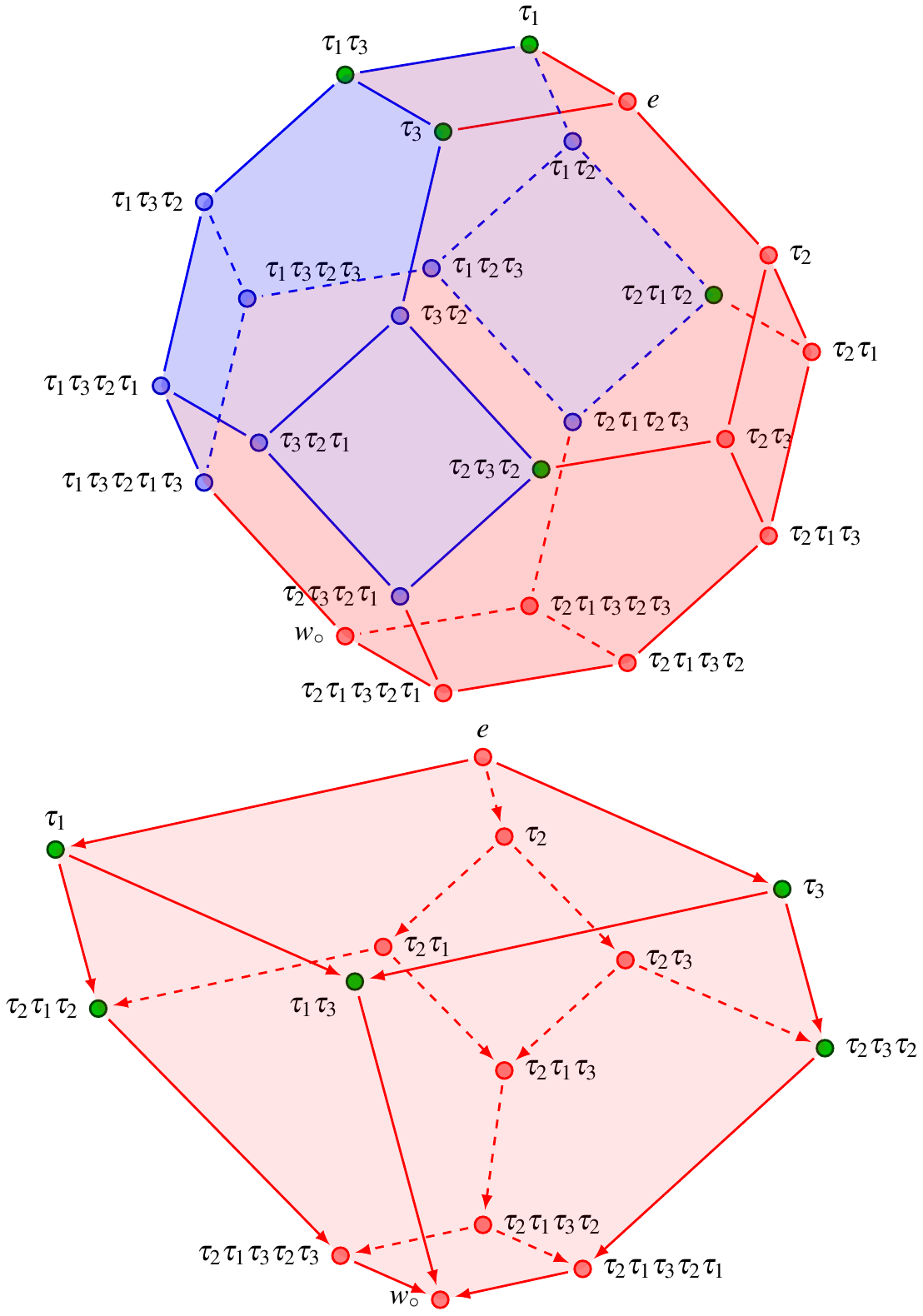}
\caption{{\bf The associahedron as realized by F.~Chapoton, S.~Fomin and A.~Zelevinsky.} On the top is the classical permutahedron for the symmetric group $S_4$, see Example~\ref{H_ex:A3Coxeter},  and on the bottom is associahedron $\textnormal{Asso}^{\pmb a}_{\tau_2\tau_1\tau_3}(S_4)$ which corresponds to the associahedron realized by Chapoton, Fomin and Zelevinsky.  The facets of $\textnormal{Asso}^{\pmb a}_{\tau_2\tau_1\tau_3}(S_4)$ containing Coxeter singletons, as well as the  Coxeter singletons, are in red. In green are the vertices corresponding to the $\tau_1\tau_2\tau_3$-sortable elements. The orientations on the edges represent the $\tau_2\tau_1\tau_3$-Cambrian lattice.}
\label{H_fig:AssA3-2}
\end{center}
\end{figure}

\begin{figure}
\begin{center}
  \includegraphics[scale=1.]{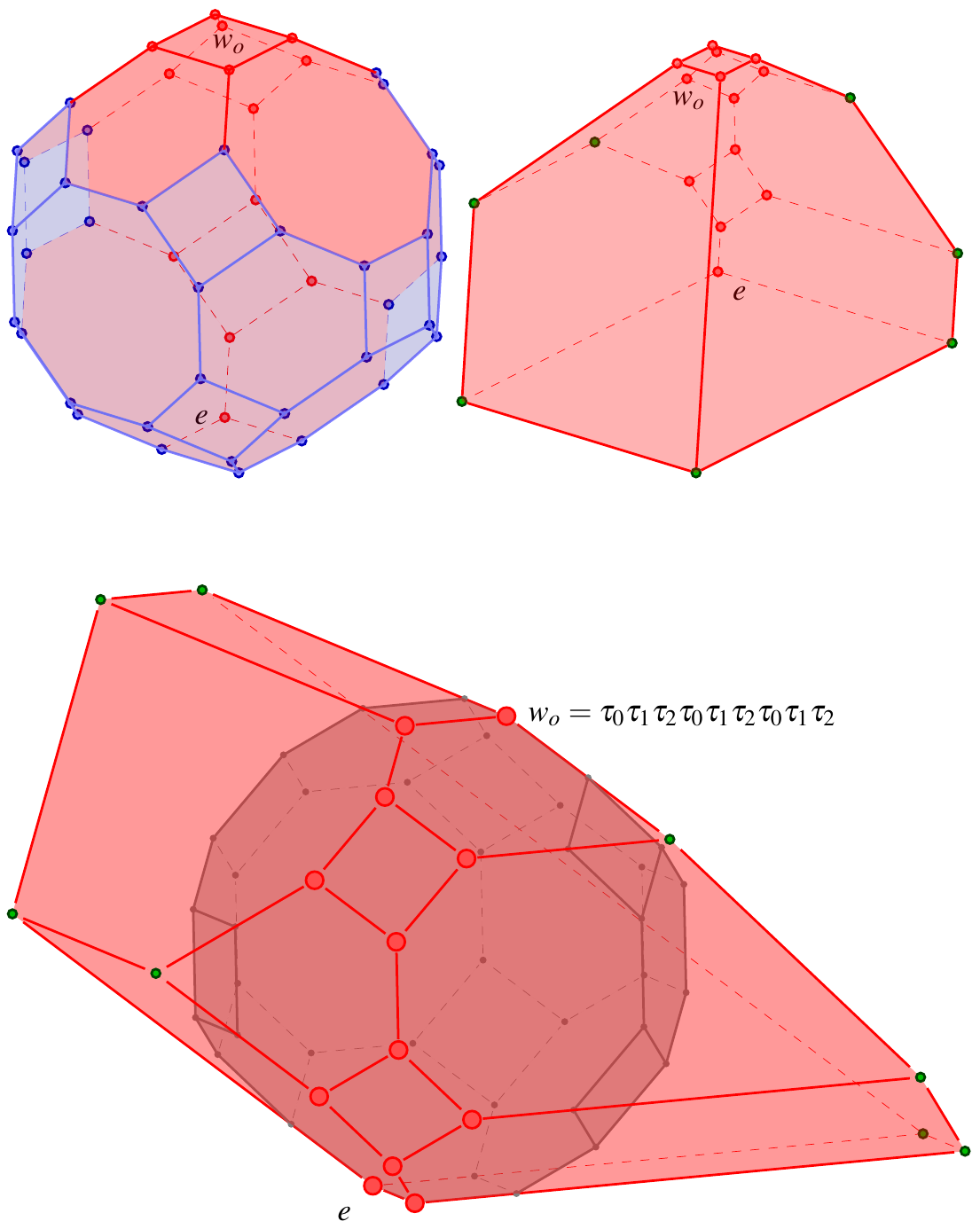}
\caption{{\bf The cyclohedron.} Two Coxeter generalized associahedra for the hyperoctahedral $W'_3$ (see Example~\ref{H_ex:BnCoxeter}). On the top left is $\textnormal{Perm}^{\pmb a}(W'_3)$, with $w_o=\tau_1\tau_2s_0\tau_1\tau_2s_0\tau_1\tau_2s_0$, and on the top right is $\textnormal{Asso}^{\pmb a}_{\tau_1\tau_2s_0}(W'_3)$. On the bottom we can see through a ghostly $ \textnormal{Asso}^{\pmb a}_{s_0\tau_1\tau_2}(W'_3)$ the permutahedron $\textnormal{Perm}^{\pmb a}(W'_3)$ in gray. In both examples, the facets containing Coxeter singletons, as well as the  Coxeter singletons, are in red.}
\label{H_fig:AssB3}
\end{center}
\end{figure}

\begin{figure}
\begin{center}
  \includegraphics[scale=1.]{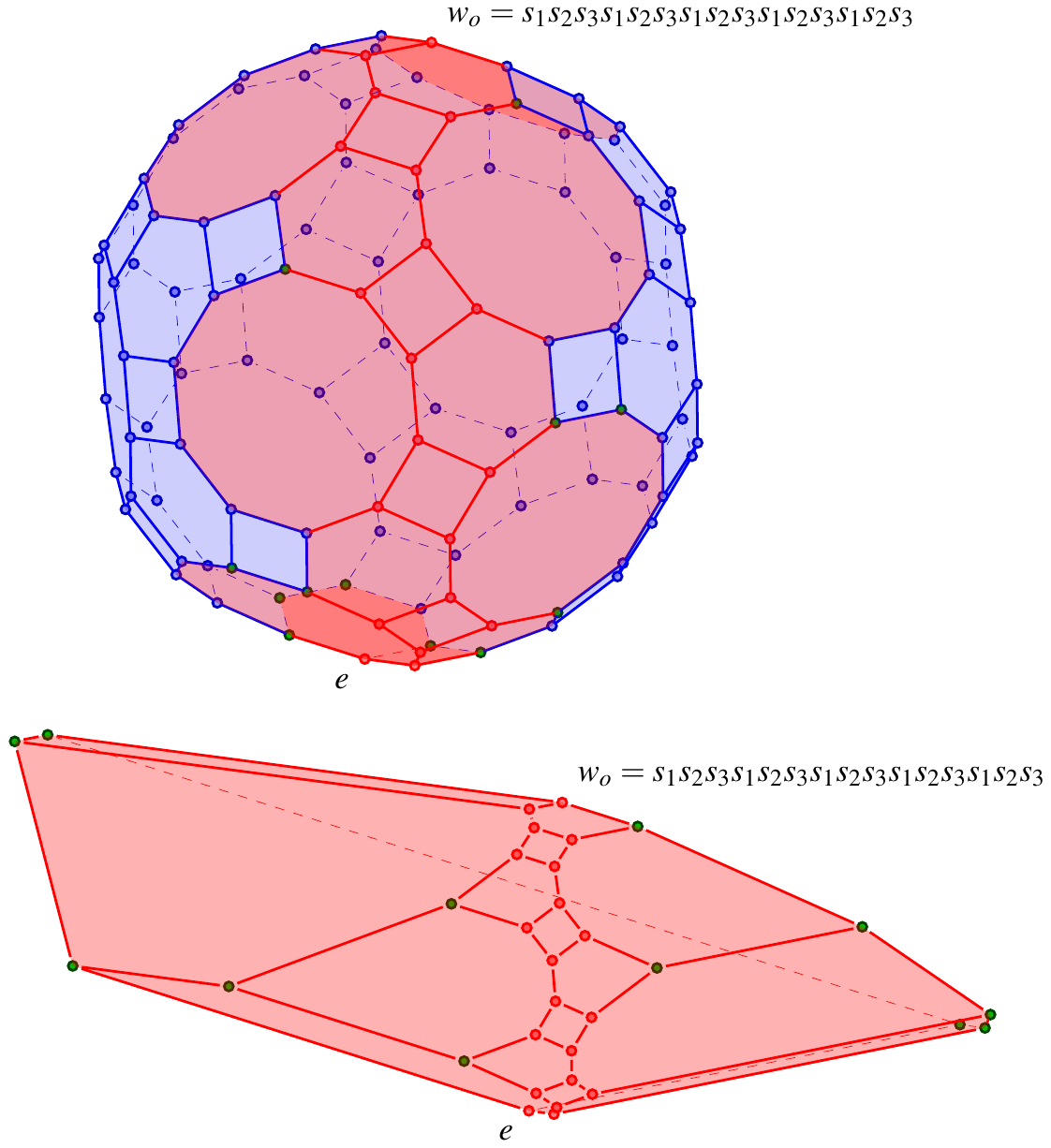}
\caption{On the top is a $W$-permutahedron for the isometry group $W(H_3)$ of the dodecahedron, see Example~\ref{H_ex:H3Coxeter},  and on the bottom is a Coxeter generalized associahedra $\textnormal{Asso}^{\pmb a}_{s_1s_2s_3}(W(H_3))$. The facets containing Coxeter singletons, as well as the  Coxeter singletons, are in red. Note that the Coxeter generalized associahedron is reduced to half of its original size to fit in the picture; the size of the $W$-permutahedron is unchanged. }
\label{H_fig:AssH3-1}
\end{center}
\end{figure}

\begin{figure}
\begin{center}
  \includegraphics[scale=1.]{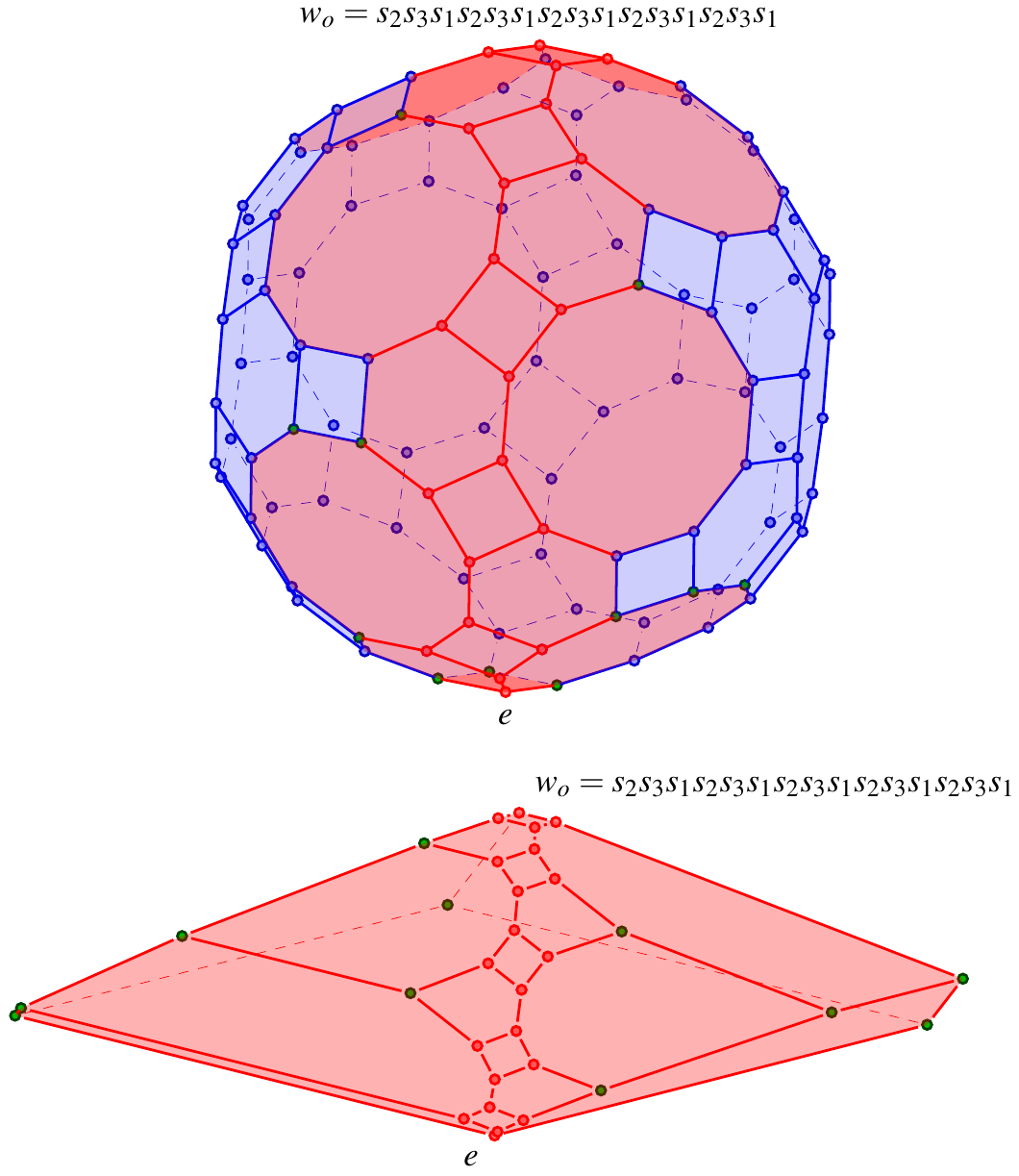}
\caption{On the top is a $W$-permutahedron for the isometry group $W(H_3)$ of the dodecahedron, see Example~\ref{H_ex:H3Coxeter},  and on the bottom is a Coxeter generalized associahedra $\textnormal{Asso}^{\pmb a}_{s_2s_3s_1}(W(H_3))$. The facets containing Coxeter singletons, as well as the  Coxeter singletons, are in red. Note that the Coxeter generalized associahedron is reduced to half of its original size to fit in the picture; the size of the $W$-permutahedron is unchanged. }
\label{H_fig:AssH3-2}
\end{center}
\end{figure}

A natural question we can ask is how many of the realizations are isometric? For instance, we observe that the Coxeter generalized associahedra shown in Figures~\ref{H_fig:AssA3-1},~\ref{H_fig:AssA3-2},~\ref{H_fig:AssB3},~\ref{H_fig:AssH3-1}~and~\ref{H_fig:AssH3-2} are not isometric. Let us briefly explain how it works.

An automorphism of the set of simple generators $S$ is a
bijection~$\mu$ on~$S$ such that the order of~$\mu(s)\mu(t)$
equals the order of $st$ for all $s,t \in S$. In particular, $\mu$
induces an automorphism on~$W$. 

\begin{proposition}[\cite{H_BeHoLaTh11}]\label{H_cor:Main} Let~$c_1$, $c_2$ be two Coxeter elements
 in $W$. Suppose that $\Hscalprod{\pmb a}{\alpha}=\Hscalprod{\pmb a}{\beta}$ for all $\alpha,\beta\in
\Delta$. Then the following statements are equivalent.
   \begin{enumerate}
      \item $\textnormal{Asso}_{c_1}(W)=\varphi\left(\textnormal{Asso}_{c_2}(W)\right)$ for some linear isometry~$\varphi$ on~$V$.
      \item There is an automorphism~$\mu$ of $S$ such that $\mu(c_2) = c_1$ or~$\mu(c_2)=c_1^{-1}$.
   \end{enumerate}
\end{proposition}

A more general statement, that is without the conditions $\Hscalprod{\pmb a}{\alpha}=\Hscalprod{\pmb a}{\beta}$, can be found in \cite{H_BeHoLaTh11}.

 \subsection{Integer coordinates}

An important subclass of finite reflection groups are {\em Weyl groups} which are linked with the theory of semi-simple Lie algebras. A finite reflection group $W$ is a {\em Weyl group} if it stabilizes a lattice in $V$, that is, a $\mathbb Z$-span of a basis of $V$.  For any Weyl group $W$, there are particular choices of root systems which are called {\em crystallographic}: a root system~$\Phi$
for~$W$ is crystallographic if for any two roots~$\alpha,\beta\in\Phi$ we have $s_\alpha(\beta)=\beta + \lambda\alpha$ for some~$\lambda\in\mathbb Z$. In that case, the simple roots~$\Delta$ span a $\mathbb Z$-lattice~$L$.  (For more details see \cite[\S2.8 and \S2.9]{H_Hu90}). Note that not all the root systems for Weyl groups are crystallographic. 

\begin{proposition}[\cite{H_HoLaTh11}]\label{H_thm:IntegerCoordinates}
   Let~$\Phi$ be a crystallographic root system for the Weyl group~$W$ and~$c$ be a Coxeter element of~$W$. Suppose that  $\pmb a\in V_0=\textnormal{span}\,(\Delta)$ has integer coordinates in $\Delta$. Then the vertices of $\textnormal{Perm}^{\pmb a}(W)$ and of $\textnormal{Asso}^{\pmb a}_c(W)$  have integer coordinates.
\end{proposition}

\begin{problem} If $\pmb a$ is not contained in $V_0$ but has integer coordinates in a $\mathbb Z$-lattice of $V$, show that $\textnormal{Perm}^{\pmb a}(W)$ and of $\textnormal{Asso}^{\pmb a}_c(W)$  have integer coordinates.
\end{problem}

This problem is supported by the following example.

\begin{example}
The classical permutahedron $\Pi_{n-1}$ obtained from the reflection group $S_n$ was presented in Example~\ref{H_ex:SnPerm}:  our construction applies to $\Pi_{n-1}$, see Remark~\ref{H_ex:HoLa}. In this setting $\pmb a=(1,2,\dots,n)\in V_{\pmb a}$ has integer coordinates. The straightforward idea to apply the last theorem to this setting is to look at the orthogonal projection onto $V_0$, which unfortunately  is $(1-(n+1)/2,2-(n+1)/2,\dots, n-(n+1)/2)$ and does not have integer coordinates. However, $S_n$ fixes the lattice spanned by the   canonical basis of $\mathbb R^n$, and it was shown in \cite{H_HoLa07} that the coordinates of $c$-generalized associahedra obtained from $\Pi_{n-1}$ are integers. We sill cannot explain this phenomenon.

\begin{problem} Let $W$ be a Weyl group and $L$ be a lattice stable under the action of $\Phi$. Suppose that  $\pmb a\in V$ has integer coordinates in $L$. Do $\textnormal{Perm}^{\pmb a}(W)$ and $\textnormal{Asso}^{\pmb a}_c(W)$  have integer coordinates?
\end{problem}
\end{example}

\section{Further developments}

We presented in this text many open problems. However, possible developments should be mentioned. First, Coxeter generalized associahedra arises in the theory of finite cluster algebra theory, so it is natural to ask the following question.

\begin{problem} Is it possible to build finite type cluster algebras from Coxeter generalized associahedra?
\end{problem}

Also, natural questions arising from  the theory of polytopes are still open, questions such as the following.

\begin{problem}
\label{H_pr:pr10}
 Compute the volume of  permutahedra and of Coxeter generalized associahedra for any finite reflection group.
\end{problem}

\begin{problem}
\label{H_pr:pr11}
 If $W$ is a Weyl group and $\Phi$ is crystallographic, compute the number of points with integers coordinates contained in $\textnormal{Asso}_c^{\pmb a}(W)$ and $\textnormal{Perm}^{\pmb a}(W)$.
\end{problem}

\index{permutahedron!generalized permutahedra}
A.~Postnikov introduced in \cite{H_Po09} the concept of {\em generalized permutahedra}, which provide tools to start answering Problem~\ref{H_pr:pr10} and Problem~\ref{H_pr:pr11}.  Generalized permutahedra are polytopes obtained from the classical permutahedron $\Pi_{n-1}$ by `nicely moving some facets', see \cite{H_PoReWi08} for details. In our setting, we are obtaining Coxeter generalized associahedra by `{\em re}moving some facets'. One of the interest of A.~Postnikov's approach is to express the classical permutahedron, as well as the associahedron realized by J.-L.~Loday and S.~Shnider-S.~Sternberg, as the {\em Minkowski sum} of faces of standard simplicies. Recently, C.~Lange expressed generalized associahedra obtained from $S_n$ as Minkowski sums {\em and differences} of faces of standard simplicies~\cite{H_La11}. This result suggests extending the framework of generalized permutahedra to include subtractions of faces of standard simplicies.  These questions were never explored in the context of finite reflection groups; one of the reasons is that we have, at this point, no idea what are the `good' objects to consider in place of the standard simplex. 

\begin{problem} For a finite reflection group $W$, find a suitable framework to define {\em Coxeter generalized permutahedra} and then express Coxeter generalized associahedra as Minkowski sums in this framework. 
\end{problem}

\index{associahedron!multi-associahedra}
Finally, the work of C.~Ceballos and J.-P.~Labb\'e and C.~Stump~\cite{H_CeLaSt11} and V.~Pilaud and C.~Stump~\cite{H_PiSt11} on {\em spherical subword complexes} suggest that there is a more general family of polytopes, arising together with generalizations of the Tamari lattice, which include Coxeter generalized associahedra: they are called {\em generalized multi-associahedra}. This family is still awaiting  polytopal realizations.

\begin{problem}[\cite{H_CeLaSt11}] Give a polytopal realization of generalized multi-associahedra.
\end{problem}

%%%%%%%%
\begin{acknowledgement}
The author is more than grateful to Carsten Lange for allowing him to use some of the pictures he made for the articles \cite{H_HoLaTh11}, and also to Jean-Philippe Labb\'e for providing the original picture in TikZ which are the bases for the pictures of permutahedra and associahedra of the symmetric group $S_4$, the hyperoctahedral group $W'_3$ and the group $W(H_3)$ presented in this article.

The author wishes also to thank Jean-Philippe Labb\'e, Carsten Lange, Vincent Pilaud and Christian Stump for many very interesting conversations which mostly took place at LaCIM (Laboratoire de Combinatoire et d'Informatique Math\'ematique) in Montr\'eal during the summer of 2011. The author thanks Franco Saliola for his comments on a preliminary version of this text.

And finally, the author expresses his deepest gratitude to Jean-Louis Loday to have asked him the question {\em how to realize the cyclohedron from a permutahedron of the hyperoctahedral group} while he was working on his PhD at the {\em Universit\'e de Strasbourg}. 
\end{acknowledgement}
%%%%%%%%

%\input{referenc}

%\bibliographystyle{tamari}
%\bibliography{MasterBibliography}

% After completion of the final version you copy the contents of the file "author.bbl" between

\providecommand{\bysame}{\leavevmode\hbox to3em{\hrulefill}\thinspace}
\providecommand{\MR}{\relax\ifhmode\unskip\space\fi MR }
\providecommand{\MRhref}[2]{%
  \href{http://www.ams.org/mathscinet-getitem?mr=#1}{#2}
}
\providecommand{\href}[2]{#2}

\index{associahedron!classical associahedron}
\index{Tamari lattice}

\index{associahedron!cyclohedron}

\end{document}